\documentclass{article}%
\usepackage{amssymb}
\usepackage{amsmath}
\usepackage{amssymb}
\usepackage{amsmath}
\usepackage{epsfig}
\usepackage{graphicx}
\usepackage{color}
\usepackage[utf8]{inputenc}
\usepackage[english]{babel}
\usepackage{caption}
\usepackage{graphicx}
\usepackage{placeins}
\usepackage{float}
\usepackage{amsfonts}%
\providecommand{\U}[1]{\protect\rule{.1in}{.1in}}

\begin{document}

\begin{center}
\textbf{A Powerful Robust Cubic Hermite Collocation Method for the Numerical
Calculations and Simulations of the Equal Width Wave Equation}

\textbf{\textrm{Sel\c{c}uk KUTLUAY}}$^{\mathbf{\mathrm{a}}}$\textbf{\textrm{,
Nuri Murat YA\u{G}MURLU}}$^{\mathbf{\mathrm{a}}}$,\textbf{\textrm{{\ }Ali
Sercan KARAKA\c{S}}}$^{\mathrm{a}}$\textbf{\textrm{{\textrm{\newline}}}%
}$^{\mathbf{\mathrm{{\mathrm{a}}}}}$%
\textbf{\textrm{{\textbf{\textrm{\textbf{In\"{o}n\"{u} University, }%
}Department of Mathematics, Malatya, 44280, TURKEY.\newline}}}}

\textbf{\textrm{{e-mail:
\textbf{\textrm{selcuk.kutluay@\textbf{\textrm{inonu.edu.tr}}}}}}
\textrm{{ORCID: }}}https://orcid.org/0000-0001-9610-504X

\textbf{\textrm{{e-mail:
\textbf{\textrm{murat.yagmurlu@\textbf{\textrm{inonu.edu.tr}}}} \ ORCID: }}}https://orcid.org/0000-0003-1593-0254

\textbf{\textrm{e-mail: ali\_sercan\_44@hotmail.com ORCID:}}https://orcid.org/0000-0001-8622-1127
\end{center}

\section{Abstract}

\ In this article, non-linear Equal Width-Wave (EW) equation will be
numerically solved . For this aim, the non-linear term in the equation is
firstly linearized by Rubin-Graves type approach. After that, to reduce the
equation into a solvable discretized linear algebraic equation system which is
the essential part of this study, the Crank-Nicolson type approximation and
cubic Hermite collocation method are respectively applied to obtain the
integration in the temporal and spatial domain directions. To be able to
illustrate the validity and accuracy of the proposed method, six test model
problems that is single solitary wave, the interaction of two solitary waves,
the interaction of three solitary waves, the Maxwellian initial condition,
undular bore and finally soliton collision will be taken into consideration
and solved. Since only the single solitary wave has an analytical solution
among these solitary waves, the error norms $\mathit{L}_{\infty}$ and $L_{2}$
are computed and compared to a few of the previous works available in the
literature. Furthermore, the widely used three invariants $I_{1},$ $I_{2}$ and
$I_{3}$ of the proposed problems during the simulations are computed and
presented. Beside those, the relative changes in those invariants are
presented. Also, a comparison of the error norms $\mathit{L}_{\infty}$ and
$L_{2}$ and these invariants obviously shows that the proposed scheme produces
better and compatible results than most of the previous works using the same
parameters. Finally, von Neumann analysis has shown that the present scheme is
unconditionally stable.

\textbf{Keywords:} Equal width-wave equation, cubic hermite collocation
method, solitary waves, stability analysis, Crank-Nicolson type approximation,
Rubin-Graves type linearization.

\textbf{AMS classification:} {65L60, 65N35, 74J35, 65D07.}

\section{\textbf{Introduction}}

Scientists encounter many physical phenomena occurring in nature and they
generally express those phenomena by algebraic, differential or integral
equations. Non-linear evolution equations are such a commonly and widely
utilized around us in order to describe complex phenomena in various areas of
sciences, however they are taken for granted. When those types of phenomena
are investigated in detail, it is seen that most of the nonlinear phenomena
which have a crucial role in science and mathematics are generally modeled by
non-linear partial differential equations (PDEs). In general, it is difficult
and troublesome to investigate and find exact solutions of initial and
boundary value problems consisting of non-linear PDEs. Actually,
scientists\ agree that there is no such a method, scheme or technique yet, it
is necessary to deal with almost every type of those equations in itself and
solve it. Because of this reason, numerical solutions are usually preferred
instead of their exact ones. Thus, many researchers are concentrated on
approximate methods and techniques to obtain numerical solutions of non-linear
PDEs. One of such equations is widely known as EW equation. This equation is
usually seen as an alternative way of defining of Korteweg-de Vries (KdV)
equation. The EW equation was firstly proposed and derived by Morrison
\textit{et al.}\cite{bl01} and is utilized as an alternative way of defining
KdV equation and presented in the following form
\begin{equation}
U_{t}+UU_{x}-\mu U_{xxt}=0, \label{1}%
\end{equation}
where $\mu$ stands for a positive parameter and the subscripts $x$ and $t$
denote partial differentiation with respect to spatial and temporal
dimensions, respectively.

There have been several analytical and numerical works about the EW equation
which has solutions showing soliton like solutions and illustrates an
equilibrium condition between nonlinear and dispersive effects available
inheritenly in the nature of the phenomena. In recent years, several studies
as in Refs \cite{bl02,bl03,bl04,bl05} can be found in the literature for the
analytical solutions of the equation. Whereas, several scholars seek numerical
solutions of the EW equation. For example, Ya\u{g}murlu and Karaka\c{s}
\cite{ek1} have found {\normalsize approximate solutions of the EW equation
using cubic trigonometric collocation finite element method based on
Rubin-Graves type linearization. }Among others, Haar wavelet
method\cite{bl010}, collocation method\cite{bl09}, Petrov-Galerkin
method\cite{bl3gardnerayoup}, least-squares method\cite{bl3zaki}, radial basis
function based pseudo-spectral method\cite{bl013}, linearized implicit
finite-difference method\cite{bl014}, lumped Galerkin method\cite{bl3esen},
explicit finite difference methods\cite{ek2}, multi-quadric
quasi-interpolation method\cite{ek5} and fully implicit finite difference
method\cite{bl016} are applied to get approximate solutions of the EW equation.

The proposed method to be used in this study is a mixture of the orthogonal
collocation method and the finite element method, where the cubic hermite
polynomials are used as a trial function. Since these polynomials satisfy the
continuity conditions for a trial function and its first and second order
derivatives at nodal points, they produce solutions with continuous
derivatives throughout the domain of the problem.

In this method, the solution region is firstly split into a number of
elements, and next orthogonal collocation is used in each one of these
elements, setting the residue equal to zero at two interior nodal points.
Nodal points have a key role in the discretization process of the equation
with respect to $x$. For the present method, the roots of orthogonal
polynomials such as the second degree Legendre and Chebyshev polynomials are
usually taken as collocation points. Arora et al.\cite{arora2007} have used
the roots of Legendre polynomials at interior collocation points and
illustrated that those polynomials present results having less error than
Chebyshev polynomials. In addition, they observed that while Chebyshev
polynomials produce better results only at cups, Legendre polynomials produce
better results on the average as well as at the cups.

In this work, we will present numerical solutions and simulations of the EW
equation using cubic Hermite B-spline collocation method with the help of
Crank-Nicolson type approximation. Truly, collocation method based on various
B-splines is generally utilized to find approximate solutions of non-linear
PDEs. Several scholars have utilized the collocation method based on various
base functions such as classical B-splines, exponential and radial base
functions and trigonometric B-splines. Regarding the article itself and its
details, one can refer to the articles
\cite{2013ganaie,2013mittal,ganaiekukreja,2016ganaie,2018kaur,2018arora,2020arora,2020yousaf,2021kaur,2021nikolay,2021kumari,2021kumaria}
and the references in it.

The present paper has been divided into seven sections. The first one is an
introduction to the Cubic Hermite Collocation Method (CHCM). A brief
description of the EW equation is given in Section two. Sections three and
four detail the implementation of the proposed scheme. Section five is about
the stability analysis of the scheme. Section six includes comparatively the
numerical results and simulations obtained by solving six test problems using
the present method. The last section, which is Section seven, is dedicated to
a brief conclusion with a future work.

\section{\textbf{Implementation of the method for spatial discretization}}

In this article, the EW equation in the following form is considered
\[
U_{t}+UU_{x}-\mu U_{xxt}%
=0,\ \ \ \ \ \ \ \ \ \ \ \ \ \ \ \ \ \ \ \ \ \ \ \ \ \ -\infty<x<\infty
\]
having the physical boundary conditions $U\rightarrow0$ when $x\rightarrow
\pm\infty,$ where $x$ is the spatial, $t$ is temporal coordinate and $\mu$ is
a positive constant. During the numerical computations of the problems to be
considered in Numerical Examples Section, the suitable boundary conditions are
going to be taken as%

\begin{equation}%
\begin{array}
[c]{c}%
U(a,t)=0,~~~~~~U(b,t)=0,\\
U_{x}(a,t)=0,~~~~~~U_{x}(b,t)=0.
\end{array}
\label{bc}%
\end{equation}
\qquad To apply a numerical method, as in general, consider that spatial
domain is chosen as the finite interval $[a,b]$ \ and then is split into $N$
finite elements with equal lengths at the collocation points $x_{i}%
$,~~$i=0(1)N$ such that $a=x_{0}<x_{1}\cdots<x_{N}=b$ and $h=x_{i+1}-x_{i}$.
The cubic hermite base functions $H_{j}$ $(j=1(1)N+1)$ are given as
\cite{ganaiekukreja}
\begin{equation}
H_{2j-1}\left(  x\right)  =\frac{1}{h^{3}}\left\{
\begin{array}
[c]{r}%
\left(  x-x_{j-1}\right)  ^{2}\left[  3h-2\left(  x-x_{j-1}\right)  \right]
,\qquad x_{j-1}\leq x\leq x_{j}\\
\left[  h-\left(  x-x_{j}\right)  \right]  ^{2}\left[  h-2\left(
x-x_{j}\right)  \right]  ,\qquad x_{j}\leq x\leq x_{j+1}\\
0,\qquad\qquad\qquad\text{otherwise}%
\end{array}
\right.  \label{hermit1}%
\end{equation}

\begin{equation}
H_{2j}\left(  x\right)  =\frac{1}{h^{3}}\left\{
\begin{array}
[c]{r}%
-h\left(  x-x_{j-1}\right)  ^{2}\left[  h-\left(  x-x_{j-1}\right)  \right]
,\qquad x_{j-1}\leq x\leq x_{j}\\
h\left(  x-x_{j}\right)  \left[  h-\left(  x-x_{j}\right)  \right]
^{2},\qquad x_{j}\leq x\leq x_{j+1}\\
0,\qquad\qquad\qquad\text{otherwise.}%
\end{array}
\right.  \label{hermit2}%
\end{equation}

An approximation solution $U_{N}(x,t)$ to $U(x,t)$ is written by means of the
cubic hermite B-splines%

\begin{equation}
U_{N}(x,t)\approx U(x,t)\approx%
{\displaystyle\sum\limits_{j=1}^{N}}
a_{j+2k-2}\left(  t\right)  H_{ji} \label{4}%
\end{equation}
in which $a$'s are time dependent parameters to be found, $k$ is the number of
elements and $i=1,2$. If the second order Legendre quadrature points\ $\eta
_{ji}\ $are chosen for each subinterval $\left[  x_{j},x_{j+1}\right]  $ ,
then the Legendre quadrature points are taken as follows%

\begin{equation}
\eta_{ji}=\frac{x_{j-1}+x_{j}}{2}+\left(  -1\right)  ^{i}\frac{h_{j}}%
{2\sqrt{3}},\quad2\leq j\leq N+1,\quad1\leq i\leq2. \label{GL}%
\end{equation}
When the following shifted Legendre polynomial roots are used in Eq. (\ref{GL})%

\[
\xi_{1}=\frac{1}{2}\left(  1+\frac{1}{\sqrt{3}}\right)  ,\quad\xi_{2}=\frac
{1}{2}\left(  1-\frac{1}{\sqrt{3}}\right)
\]
one gets%

\[
\frac{\eta_{j1}-x_{j}}{h_{j}}=-\xi_{1},\quad\frac{\eta_{j2}-x_{j}}{h_{j}}%
=-\xi_{2}%
\]
B{\normalsize ut, if }Chebyshev polynomial is chosen the following roots%

\[
\xi_{1}=\frac{1}{2}\left(  1+\frac{1}{\sqrt{2}}\right)  ,\quad\xi_{2}=\frac
{1}{2}\left(  1-\frac{1}{\sqrt{2}}\right)
\]
are used. Throughout the article, both Legendre and Chebyshev polynomial roots
are used for the numerical computations.

In this method, after discretization, a new coordinate variable $\xi$ is
defined in each element such that $\xi=(x-x_{k})/h$. Thus, the variable $x$
changes in the range $\left[  x_{k},x_{k+1}\right]  $, while the new variable
$\xi$ changes in the range of $\left[  0,1\right]  $. Thus using the
transformation $x=h\xi+x_{k}$, the following equations are obtained%

\begin{align*}
H_{1}\left(  \xi\right)   &  =\left(  1+2\xi\right)  \left(  1-\xi\right)
^{2},\quad H_{2}\left(  \xi\right)  =\xi\left(  1-\xi\right)  ^{2}h\qquad\\
H_{3}\left(  \xi\right)   &  =\xi^{2}\left(  3-2\xi\right)  ,\quad
\text{\ \ \ \ \ \ \ \ }H_{4}\left(  \xi\right)  =\xi^{2}\left(  \xi-1\right)
h
\end{align*}%
\begin{align*}
A_{1}\left(  \xi\right)   &  =6\xi^{2}-6\xi,\quad\text{\ \ }A_{2}\left(
\xi\right)  =\left(  1-4\xi+3\xi^{2}\right)  h\qquad\ \\
A_{3}\left(  \xi\right)   &  =6\xi-6\xi^{2},\quad A_{4}\left(  \xi\right)
=\left(  3\xi^{2}-2\xi\right)  h
\end{align*}%
\begin{align*}
B_{1}\left(  \xi\right)   &  =12\xi-6,\quad B_{2}\left(  \xi\right)  =\left(
6\xi-4\right)  h\qquad\text{\ \ \ \ \ \ \ \ \ \ \ \ \ }\\
B_{3}\left(  \xi\right)   &  =6-12\xi,\quad B_{4}\left(  \xi\right)  =\left(
6\xi-2\right)  h.
\end{align*}

Thus, the trial function over the $k^{th}$ element are defined as%

\[
U_{N}(x,t)=%
{\displaystyle\sum\limits_{j=1}^{N}}
a_{j+2k-2}\left(  t\right)  H_{ji}.
\]
The trial functions with their first and second order derivatives at the
collocation points in terms of local variable $\xi$ are defined as follows%

\begin{align*}
U_{N}\left(  \xi,t\right)   &  =%
{\displaystyle\sum\limits_{j=1}^{4}}
a_{j+2k-2}\left(  t\right)  H_{j}\left(  \xi\right) \\
&  =a_{2k-1}H_{1}\left(  \xi\right)  +a_{2k}H_{2}\left(  \xi\right)
+a_{2k+1}H_{3}\left(  \xi\right)  +a_{2k+2}H_{4}\left(  \xi\right)
\end{align*}%
\begin{align*}
U_{N}^{^{\prime}}\left(  \xi,t\right)   &  =\frac{1}{h}%
{\displaystyle\sum\limits_{j=1}^{4}}
a_{j+2k-2}\left(  t\right)  A_{j}\left(  \xi\right) \\
&  =\frac{1}{h}\left[  a_{2k-1}A_{1}\left(  \xi\right)  +a_{2k}A_{2}\left(
\xi\right)  +a_{2k+1}A_{3}\left(  \xi\right)  +a_{2k+2}A_{4}\left(
\xi\right)  \right]
\end{align*}%
\begin{align*}
U_{N}^{^{^{\prime\prime}}}\left(  \xi,t\right)   &  =\frac{1}{h^{2}}%
{\displaystyle\sum\limits_{j=1}^{4}}
a_{j+2k-2}\left(  t\right)  B_{j}\left(  \xi\right) \\
&  =\frac{1}{h^{2}}\left[  a_{2k-1}B_{1}\left(  \xi\right)  +a_{2k}%
B_{2}\left(  \xi\right)  +a_{2k+1}B_{3}\left(  \xi\right)  +a_{2k+2}%
B_{4}\left(  \xi\right)  \right]  .
\end{align*}
Here $A_{1},A_{2},A_{3},A_{4}$ ve $B_{1},B_{2},B_{3},B_{4}$ are the first and
second order derivatives of Hermite base functions, respectively. When Eqs.
(\ref{hermit1}) and (\ref{hermit2}) are used at the nodal points, the
following approximate solutions are found%

\begin{equation}%
\begin{array}
[c]{c}%
\begin{array}
[c]{c}%
\text{ \ \ \ \ }U_{i}=U_{N}\left(  \xi_{i},t\right)  =a_{2k-1}H_{1i}%
+a_{2k}H_{2i}+a_{2k+1}H_{3i}+a_{2k+2}H_{4i}\\
\text{ \ }hU_{i}^{\prime}=U_{N}^{^{\prime}}\left(  \xi_{i},t\right)
=a_{2k-1}A_{1i}+a_{2k}A_{2i}+a_{2k+1}A_{3i}+a_{2k+2}A_{4i}%
\end{array}
\\
h^{2}U_{i}^{^{\prime\prime}}=U_{N}^{^{^{\prime\prime}}}\left(  \xi
_{i},t\right)  =a_{2k-1}B_{1i}+a_{2k}B_{2i}+a_{2k+1}B_{3i}+a_{2k+2}B_{4i}%
\end{array}
\label{chcmbaz}%
\end{equation}
where$\ H_{ji}=H_{j}\,\left(  \xi_{i}\right)  $, $A_{ji}=A_{j}\,\left(
\xi_{i}\right)  $ and $B_{ji}=B_{j}\,\left(  \xi_{i}\right)  $ for $i=1,2.$

During the solution process, firstly, forward finite difference approximation
for temporal integration and then finite element collocation method using
cubic Hermite B-spline basis functions for spatial integration will be
implemented. In fact, the implementation of the presented method based on
Hermite B-splines are more efficient because of their several crucial
characteristics such as easy storage and manipulations in computers.

It is worth to note that both of linear and non-linear algebraic equations
systems found using any B-splines are usually well-conditioned and let the
required parameters be determined quite easily. Furthermore, when obtaining
the approximations by B-splines, one mostly doesn't encounter numerical
instability. Moreover, the matrix systems found by B-splines are in general
sparse band matrixes and easy to be implemented on digital computers.

\section{\textbf{Implementation of the method for temporal discretization}}

At the moment, we will discretize the EW equation $(\ref{1})$ given as
\[
U_{t}+UU_{x}-\mu U_{xxt}=0.
\]
To do so, we have firstly implemented the Crank-Nicolson type approximation to
Eq. $(\ref{1})$ to get the following discretized scheme%
\begin{equation}
\frac{U^{n+1}-U^{n}}{\Delta t}+\frac{(UU_{x})^{n}+(UU_{x})^{n+1}}{2}-\mu
\frac{\left(  U_{xx}\right)  ^{n+1}-\left(  U_{xx}\right)  ^{n}}{\Delta t}=0.
\label{15}%
\end{equation}
Then, linearizing the nonlinear term $\left(  UU_{x}\right)  ^{n+1}$ in Eq.
$(\ref{15})$ by virtue of the Rubin-Graves approximation\cite{bl017}
\begin{equation}
\left(  UU_{x}\right)  ^{n+1}=U_{x}^{n+1}U^{n}+U^{n+1}U_{x}^{n}-U_{x}^{n}%
U^{n}, \label{rb}%
\end{equation}
and substituting $(\ref{rb})$ into $(\ref{15})$, one gets the following
recursive formula to find next time level unknowns%
\begin{equation}
U^{n+1}(\frac{1}{\Delta t}+\frac{1}{2}U_{x}^{n})+\frac{1}{2}U_{x}^{n+1}%
U^{n}-\frac{\mu}{\Delta t}U_{xx}^{^{n+1}}=\frac{1}{\Delta t}U^{n}-\frac{\mu
}{\Delta t}U_{xx}^{n}. \label{ozyin}%
\end{equation}
When the cubic Hermite base functions and their derivatives given in Eq.
(\ref{chcmbaz}) are used in Eq. (\ref{ozyin}), the following iterative formula
is obtained%
\begin{equation}%
\begin{array}
[c]{l}%
\text{ \ }\left[  a_{2k-1}^{n+1}H_{1i}+a_{2k}^{n+1}H_{2i}+a_{2k+1}^{n+1}%
H_{3i}+a_{2k+2}^{n+1}H_{4i}\right]  \left[  \frac{1}{\Delta t}+\frac
{a_{2k-1}^{n}A_{1i}+a_{2k}^{n}A_{2i}+a_{2k+1}^{n}A_{3i}+a_{2k+2}^{n}A_{4i}%
}{2h}\right] \\
\\
+\left[  \frac{a_{2k-1}^{n+1}A_{1i}+a_{2k}^{n+1}A_{2i}+a_{2k+1}^{n+1}%
A_{3i}+a_{2k+2}^{n+1}A_{4i}}{h}\right]  \left[  \frac{a_{2k-1}^{n}%
H_{1i}+a_{2k}^{n}H_{2i}+a_{2k+1}^{n}H_{3i}+a_{2k+2}^{n}H_{4i}}{2}\right] \\
\\
-\frac{\mu}{\Delta t}\left[  \frac{a_{2k-1}^{n+1}B_{1i}+a_{2k}^{n+1}%
B_{2i}+a_{2k+1}^{n+1}B_{3i}+a_{2k+2}^{n+1}B_{4i}}{h^{2}}\right] \\
\\
=\left[  \frac{a_{2k-1}^{n}H_{1i}+a_{2k}^{n}H_{2i}+a_{2k+1}^{n}H_{3i}%
+a_{2k+2}^{n}H_{4i}}{\Delta t}\right]  -\frac{\mu}{\Delta t}\left[
\frac{a_{2k-1}^{n}B_{1i}+a_{2k}^{n}B_{2i}+a_{2k+1}^{n}B_{3i}+a_{2k+2}%
^{n}B_{4i}}{h^{2}}\right]
\end{array}
\label{kararlilik}%
\end{equation}
in which $T$ is being the desired final time, $\Delta t=T/M$ and
$t_{n}=n\Delta t$ $(n=1(1)M).$ From Eq. (\ref{kararlilik}), a discretized
linear algebraic system of equations is obtained. These equations are
recursive relationships for the element parameters vector $\mathbf{a}%
^{n}=(a_{1}^{n},...,a_{2N+1}^{n},a_{2N+2}^{n})$ where $t_{n}=n\Delta t,$
$n=1(1)M$ until the final time $T$. Using the boundary conditions given in
Eq.(\ref{bc}) and eliminating the parameters $a_{1}^{n},a_{2N+1}^{n}$ in Eq.
(\ref{kararlilik}) as follows: From the left boundary condition $U(x_{0}%
,t)=a_{1}^{n}H_{11}+a_{2}^{n}H_{21}+a_{3}^{n}H_{31}+a_{4}^{n}H_{41}=0,$ since
$H_{21}=H_{31}=H_{41}=0$ and $H_{11}\neq0,$ the condition $a_{1}^{n}=0$ is
obtained. Similarly from the right boundary condition $U(x_{N},t)=a_{2N-1}%
^{n}H_{12}+a_{2N}^{n}H_{22}+a_{2N+1}^{n}H_{32}+a_{2N+2}^{n}H_{42}=0,$ since
$H_{12}=H_{22}=H_{42}=0$ and $H_{32}\neq0$, the condition $a_{2N+1}^{n}=0$ is obtained.

Finally, one gets a new uniquely solvable algebraic equation system in the
following matrix form%

\begin{equation}
\mathbf{L\mathbf{\boldsymbol{a}}}^{n+1}\mathbf{\mathbf{=\boldsymbol{R}}a}%
^{n}{\normalsize .} \label{recursive}%
\end{equation}
Here the matrix $\mathbf{\mathbf{L}}$ and $\mathbf{R}$ are square $2N\times2N$
diagonal band matrices, and the matrices $\mathbf{\mathbf{\boldsymbol{a}}%
}^{n+1}$ and $\mathbf{\mathbf{\boldsymbol{a}}}^{n}$ are $2N\times1$ column matrices.

The values $\mathbf{a}_{i}$ $(i=1(1)2N)$ obtained by solving the system of
equations given by Eq.(\ref{recursive}) are found and the approximate
solutions of EW equation at the next time level are computed. This process is
repeated successively for $t_{n}=n\Delta t$ $(n=1(1)M)$ until the final time
$T$. In order to start the iterative process, the initial\ vector\ $\mathbf{a}%
^{0}$ with entries $\mathbf{a}_{i0}$ $(i=1(1)2N)$ is needed. This vector is
calculated by the initial condition presented by the governing equation.

\subsection{\textbf{The initial state}}

The initial vector $\mathbf{a}^{0}$ is found using the initial/boundary
conditions. Thus, the approximate solution in Eq. ($\ref{4})$ is written now
for the initial condition as%

\[
U(x,t)\approx U_{N}(x,t)=%
{\displaystyle\sum\limits_{j=1}^{N}}
a_{j+2k-2}^{0}\left(  t\right)  H_{ji}{\small \ }%
\]
in which the $a_{m}^{0}$'s are unknown parameters to be computed. It is
required that the initial numerical approximation $U_{N}(x,0)$ satisfify the
following conditions%

\[%
\begin{array}
[c]{c}%
U_{N}(x_{i},0)=U(x_{i},0),\ \ \ \ \ \ i=0,1,...,N\\
(U_{N})_{x}(a,0)=0,\ \ \ \ \ \ \ \ \ \ \ \ \ \ \ (U_{N})_{x}(b,0)=0.
\end{array}
\]
Thus, the matrix equation of the following form is obtained {\small
\begin{equation}
\mathbf{Wa}^{0}=\mathbf{b} \label{denkek}%
\end{equation}
} where {\small {\
\[
\mathbf{W}=\left[
\begin{array}
[c]{cccccccc}%
H_{21} & H_{31} & H_{41} &  &  &  &  & \\
H_{22} & H_{32} & H_{42} &  &  &  &  & \\
& H_{11} & H_{21} & H_{31} & H_{41} &  &  & \\
& H_{12} & H_{22} & H_{32} & H_{42} &  &  & \\
&  & \ddots & \ddots & \ddots & \ddots &  & \\
&  &  & H_{11} & H_{21} & H_{31} & H_{41} & \\
&  &  & H_{12} & H_{22} & H_{32} & H_{42} & \\
&  &  &  &  & H_{11} & H_{21} & H_{41}\\
&  &  &  &  & H_{12} & H_{22} & H_{42}%
\end{array}
\right]  ,
\]
}%
\[%
\begin{array}
[c]{l}%
\mathbf{a}{^{0}}=(a_{2},a_{3},a_{4},\ldots,a_{2N-1},a_{2N},a_{2N+2})^{T}%
\end{array}
\]
}and {\small \
\[%
\begin{array}
[c]{l}%
\mathbf{b}=(U(x_{11},0),U(x_{12},0),U(x_{21},0),\ldots,U(x_{N1},0),U(x_{N2}%
,0))^{T}.
\end{array}
\]
}

Thus the initial values required to start the numerical scheme in Eq.
(\ref{recursive}) are computed from Eq. (\ref{denkek}) and then the iteration
process is repeatedly continued until the desired final time $T$.

\section{Stability analysis}

In order to examine the stability of the linear numerical scheme
(\ref{kararlilik}), we have used von-Neumann method. For this aim,
substituting the Fourier mode%
\[
a_{j}^{n}=\xi^{n}e^{ij\varphi}%
\]
into Eq. (\ref{kararlilik}), one gets%
\begin{equation}%
\begin{array}
[c]{l}%
\xi^{n+1}e^{i(2j-1)\varphi}(\alpha_{1})+\xi^{n+1}e^{i(2j)\varphi}(\alpha
_{2})+\xi^{n+1}e^{i(2j+1)\varphi}(\alpha_{3})+\xi^{n+1}e^{i(2j+2)\varphi
}(\alpha_{4})=\\
\\
\xi^{n}e^{i(2j-1)\varphi}(\beta_{1})+\xi^{n}e^{i(2j)\varphi}(\beta_{2}%
)+\xi^{n}e^{i(2j+1)\varphi}(\beta_{3})+\xi^{n}e^{i(2j+2)\varphi}(\beta_{4})
\end{array}
\label{kararlilik3}%
\end{equation}
where $\varphi=\beta h,$ $\beta$ is the mode number, $h$ is the spatial step
size, $i=\sqrt{-1}$ and%

\begin{align*}
\alpha_{1}  &  =H_{1i}(\frac{1}{\Delta t}+\frac{a_{2k-1}^{n}A_{1i}+a_{2k}%
^{n}A_{2i}+a_{2k+1}^{n}A_{3i}+a_{2k+2}^{n}A_{4i}}{2h})\\
&  +A_{1i}(\frac{a_{2k-1}^{n}H_{1i}+a_{2k}^{n}H_{2i}+a_{2k+1}^{n}%
H_{3i}+a_{2k+2}^{n}H_{4i}}{2h})-\frac{\mu}{\Delta t}\frac{B_{1i}}{h^{2}}\\
& \\
\alpha_{2}  &  =H_{2i}(\frac{1}{\Delta t}+\frac{a_{2k-1}^{n}A_{1i}+a_{2k}%
^{n}A_{2i}+a_{2k+1}^{n}A_{3i}+a_{2k+2}^{n}A_{4i}}{2h})\\
&  +A_{2i}(\frac{a_{2k-1}^{n}H_{1i}+a_{2k}^{n}H_{2i}+a_{2k+1}^{n}%
H_{3i}+a_{2k+2}^{n}H_{4i}}{2h})-\frac{\mu}{\Delta t}\frac{B_{2i}}{h^{2}}\\
& \\
\alpha_{3}  &  =H_{3i}(\frac{1}{\Delta t}+\frac{a_{2k-1}^{n}A_{1i}+a_{2k}%
^{n}A_{2i}+a_{2k+1}^{n}A_{3i}+a_{2k+2}^{n}A_{4i}}{2h})\\
&  +A_{3i}(\frac{a_{2k-1}^{n}H_{1i}+a_{2k}^{n}H_{2i}+a_{2k+1}^{n}%
H_{3i}+a_{2k+2}^{n}H_{4i}}{2h})-\frac{\mu}{\Delta t}\frac{B_{3i}}{h^{2}}\\
& \\
\alpha_{4}  &  =H_{4i}(\frac{1}{\Delta t}+\frac{a_{2k-1}^{n}A_{1i}+a_{2k}%
^{n}A_{2i}+a_{2k+1}^{n}A_{3i}+a_{2k+2}^{n}A_{4i}}{2h})\\
&  +A_{4i}(\frac{a_{2k-1}^{n}H_{1i}+a_{2k}^{n}H_{2i}+a_{2k+1}^{n}%
H_{3i}+a_{2k+2}^{n}H_{4i}}{2h})-\frac{\mu}{\Delta t}\frac{B_{4i}}{h^{2}}%
\end{align*}

\[
\beta_{1}=\frac{H_{1i}}{\Delta t}-\frac{\mu}{\Delta t}\frac{B_{1i}}{h^{2}%
},\text{\ }\beta_{2}=\frac{H_{2i}}{\Delta t}-\frac{\mu}{\Delta t}\frac{B_{2i}%
}{h^{2}},\text{ }\beta_{3}=\frac{H_{3i}}{\Delta t}-\frac{\mu}{\Delta t}%
\frac{B_{3i}}{h^{2}},\text{ }\beta_{4}=\frac{H_{4i}}{\Delta t}-\frac{\mu
}{\Delta t}\frac{B_{4i}}{h^{2}}.
\]
Making the required algebraic manipulations in Eq.(\ref{kararlilik3}), one
obtains%
\begin{equation}
\xi=\frac{P-iQ}{R+iS} \label{ekdenk}%
\end{equation}
where%
\begin{align}
P  &  =\beta_{4}\cos2\varphi+(\beta_{1}+\beta_{3})\cos\varphi+\beta_{2},\text{
\ \ }Q=-i(-\beta_{4}\sin2\varphi+(\beta_{1}-\beta_{3})\sin\varphi)\nonumber\\
R  &  =\alpha_{4}\cos2\varphi+(\alpha_{3}+\alpha_{1})\cos\varphi+\alpha
_{2},\text{ \ \ }S=i(\alpha_{4}\sin2\varphi+(\alpha_{3}-\alpha_{1})\sin
\varphi).\nonumber
\end{align}
When the modulus of Eq. (\ref{ekdenk}) is taken, the inequality $\left\vert
\xi\right\vert \leq1$ is found, and this is the expected condition for the
numerical scheme to be unconditionally stable.

\section{Numerical experiments}

In the present section, six widely used test problems for the EW equation will
be solved and the obtained results are going to be compared to those of
existing in the literature. When the analytical solution of the test problem
exists, the validity and accuracy of the method will be controlled utilizing
the error norms $L_{2}$ and $L_{\infty}$ given as follows, respectively:
\[
L_{2}=\left(  h\sum_{i=1}^{N}\left\vert U_{i}-(U_{N})_{i}\right\vert
^{2}\right)  ^{1/2},\text{\qquad}L_{\infty}=\max_{1\leq i\leq N}\left\vert
U_{i}-(U_{N})_{i}\right\vert .
\]
\qquad In addition to these error norms, three invariants in the discrete
points, of which formulae are given as below \cite{bl018}, are computed%
\[
I_{1}=\int_{-\infty}^{\infty}Udx,\quad I_{2}=\int_{-\infty}^{\infty}\left(
U^{2}+\mu U_{x}^{2}\right)  dx,\quad I_{3}=\int_{-\infty}^{\infty}U^{3}dx.
\]
Next, the relative changes in these invariants while the program is running
are computed from
\[
I_{p}^{\ast}=\frac{I_{p}(T)-I_{p}(t_{0})}{I_{p}(t_{0})},\ p=1,2,3
\]
and also compared with their exact values. All numerical computations are made
by using both Cubic Hermite Collocation Method with Legendre roots (CHCM-L)
and Cubic Hermite Collocation Method with Chebyshev roots (CHCM-C). These
computations have been done using MATLAB R2021a on Intel (R) Core(TM) i7 8565U
CPU @1.80Ghz computer having 8 GB of RAM.

\subsection{Single solitary wave}

The first experimental problem is known as single solitary wave and it has got
an exact solution in the following form \cite{bl01}%

\begin{equation}
U(x,t)=3c\sec\text{h}^{2}\left[  k\left(  x-x_{0}-vt\right)  \right]
\label{27}%
\end{equation}
where $k=1/\sqrt{4\mu}$ is the width of the solitary wave, $\mu=1,v=c$ stands
for the velocity of the wave and $3c$ is taken as the amplitude of the wave.

Using the\ solution domain of the problem as $(x,t)\in\lbrack a,b]\times
\lbrack0,T]$, the initial condition is taken from Eq. $(\ref{27})$ at time
$t=0$ of the following form%

\[
U(x,0)=3c\sec\text{h}^{2}\left[  k\left(  x-x_{0}\right)  \right]
\]
and the boundary conditions are given by Eq. (\ref{bc}).

The exact values of the those invariants are calculated as
follows\cite{bl3gardnerayoup}%

\[
I_{1}=6\frac{c}{k},\quad I_{2}=12\frac{c^{2}}{k}+\frac{48}{5}kc^{2}\mu,\quad
I_{3}=\frac{144}{5}\frac{c^{3}}{k}.
\]
\begin{table}[ptb]
\caption{Comparison of $\ $the calculated invariants and error norms of
Problem 1 for $h=0.03$ and $k=0.05$ ($\mu=1,$ $3c=0.3,$ $x_{0}=10,$ $0\leq
x\leq30,$ $0\leq t\leq80$).}%
\label{p1t1}
{\scriptsize
\begin{tabular}
[c]{ccccccc}\hline
Method & $t$ & $I_{1}$ & $I_{2}$ & $I_{3}$ & $L_{2}\times10^{3}$ & $L_{\infty
}\times10^{3}$\\\hline
CHCM-L & $0$ & \multicolumn{1}{l}{$1.1999445724$} &
\multicolumn{1}{l}{$0.2880000252$} & $0.0576000000$ & $0.000679$ &
$0.003911$\\
& $10$ & \multicolumn{1}{l}{$1.2000134450$} &
\multicolumn{1}{l}{$0.2880000287$} & $0.0576000016$ & $0.023823$ &
$0.032148$\\
& $20$ & $1.2000387691$ & $0.2880000300$ & $0.0576000018$ & $0.032656$ &
$0.044079$\\
& $30$ & $1.2000480346$ & $0.2880000307$ & $0.0576000018$ & $0.035919$ &
$0.048468$\\
& $40$ & $1.2000512988$ & $0.2880000310$ & $0.0576000018$ & $0.037137$ &
$0.050083$\\
& $50$ & $1.2000521015$ & $0.2880000310$ & $0.0576000018$ & $0.037608$ &
$0.050678$\\
& $60$ & $1.2000513090$ & $0.2880000310$ & $0.0576000018$ & $0.037814$ &
$0.050897$\\
& $70$ & $1.2000480555$ & $0.2880000310$ & $0.0576000018$ & $0.037960$ &
$0.050978$\\
& $80$ & $1.2000388017$ & $0.2880000310$ & $0.0576000018$ & $0.038334$ &
$0.051008$\\
CHCM-C & $80$ & $1.2000388189$ & $0.2880000287$ & $0.0576000011$ & $0.040416$
& $0.051117$\\
\cite{ek1} & $80$ & \multicolumn{1}{l}{$1.1999851019$} &
\multicolumn{1}{l}{$0.2879999949$} & $0.0575999982$ & $0.024562$ &
$0.009604$\\
\cite{bl3zaki} & $80$ & \multicolumn{1}{l}{$1.1964$} &
\multicolumn{1}{l}{$0.2858$} & \multicolumn{1}{l}{$0.0569$} &
\multicolumn{1}{l}{$7.444$} & \multicolumn{1}{l}{$4.373$}\\
\cite{bl3gardnerayoup} & $80$ & \multicolumn{1}{l}{$1.1910$} &
\multicolumn{1}{l}{$0.28550$} & \multicolumn{1}{l}{$0.05582$} &
\multicolumn{1}{l}{$3.849$} & \multicolumn{1}{l}{$2.646$}\\
\cite{bl014} & $80$ & \multicolumn{1}{l}{$1.20004$} &
\multicolumn{1}{l}{$0.28799$} & \multicolumn{1}{l}{$0.0576$} &
\multicolumn{1}{l}{$0.125$} & \multicolumn{1}{l}{$0.073$}\\
\cite{bl3esen} & $80$ & \multicolumn{1}{l}{$1.19995$} &
\multicolumn{1}{l}{$0.28798$} & \multicolumn{1}{l}{$0.05759$} &
\multicolumn{1}{l}{$0.029$} & \multicolumn{1}{l}{$0.021$}\\
\cite{bl1idrisdag} & $80$ & \multicolumn{1}{l}{$1.19998$} &
\multicolumn{1}{l}{$0.28798$} & \multicolumn{1}{l}{$0.05759$} &
\multicolumn{1}{l}{$0.056$} & \multicolumn{1}{l}{$0.053$}\\
\cite{bl3dogan} & $80$ & \multicolumn{1}{l}{$1.23387$} &
\multicolumn{1}{l}{$0.29915$} & \multicolumn{1}{l}{$0.06097$} &
\multicolumn{1}{l}{$24.697$} & \multicolumn{1}{l}{$16.425$}\\
\cite{bl3bsaka} & $80$ & \multicolumn{1}{l}{$1.19999$} &
\multicolumn{1}{l}{$0.28801$} & \multicolumn{1}{l}{$0.05760$} &
\multicolumn{1}{l}{$0.003064$} & \multicolumn{1}{l}{$0.001704$}\\
\cite{bl3roshan} & $80$ & \multicolumn{1}{l}{$1.20004$} &
\multicolumn{1}{l}{$0.2880$} & \multicolumn{1}{l}{$0.0576$} &
\multicolumn{1}{l}{$0.03882$} & \multicolumn{1}{l}{$0.05151$}\\
\cite{bl3fazal} & $80$ & \multicolumn{1}{l}{$1.20004$} &
\multicolumn{1}{l}{$0.2880$} & \multicolumn{1}{l}{$0.0576$} &
\multicolumn{1}{l}{$0.03962$} & \multicolumn{1}{l}{$0.05446$}\\\hline
Analytical &  & \multicolumn{1}{l}{$1.2$} & \multicolumn{1}{l}{$0.288$} &
\multicolumn{1}{l}{$0.0576$} &  & \\\hline
\end{tabular}
\ \ \ \ \ \ \ \ \ \ }\end{table}\begin{table}[ptb]
\caption{Comparison of $\ $the calculated invariants and error norms of
Problem 1 for $h=k=0.05$ ($\mu=1,$ $3c=0.03,$ $x_{0}=10,$ $0\leq x\leq30,$
$0\leq t\leq80$).}%
\label{p1t2}
{\scriptsize
\begin{tabular}
[c]{ccccccc}\hline
Method & $t$ & $I_{1}$ & $I_{2}$ & $I_{3}$ & $L_{2}\times10^{3}$ & $L_{\infty
}\times10^{3}$\\\hline
CHCM-L & $0$ & \multicolumn{1}{l}{$0.1199943936$} &
\multicolumn{1}{l}{$0.0028800002$} & $0.0000576000$ & $0.000088$ &
$0.000391$\\
& $10$ & \multicolumn{1}{l}{$0.1199954305$} &
\multicolumn{1}{l}{$0.0028800002$} & $0.0000576000$ & $0.000339$ &
$0.000444$\\
& $20$ & $0.1199963688$ & $0.0028800002$ & $0.0000576000$ & $0.000657$ &
$0.000865$\\
& $30$ & $0.1199972178$ & $0.0028800002$ & $0.0000576000$ & $0.000948$ &
$0.001249$\\
& $40$ & $0.1199979861$ & $0.0028800002$ & $0.0000576000$ & $0.001212$ &
$0.001597$\\
& $50$ & $0.1199986812$ & $0.0028800002$ & $0.0000576000$ & $0.001451$ &
$0.001911$\\
& $60$ & $0.1199993102$ & $0.0028800002$ & $0.0000576000$ & $0.001668$ &
$0.002196$\\
& $70$ & $0.1199998793$ & $0.0028800002$ & $0.0000576000$ & $0.001864$ &
$0.002453$\\
& $80$ & $0.1200003943$ & $0.0028800002$ & $0.0000576000$ & $0.002041$ &
$0.002686$\\
CHCM-C & $80$ & $0.1200003983$ & $0.0028800002$ & $0.0000576000$ & $0.002130$
& $0.002697$\\
\cite{bl3esen} & $80$ & \multicolumn{1}{l}{$0.12000$} &
\multicolumn{1}{l}{$0.00288$} & \multicolumn{1}{l}{$0.000058$} &
\multicolumn{1}{l}{$0.003$} & \multicolumn{1}{l}{$0.002$}\\
\cite{bl3dogan} & $80$ & \multicolumn{1}{l}{$0.12088$} &
\multicolumn{1}{l}{$0.00291$} & \multicolumn{1}{l}{$0.000059$} &
\multicolumn{1}{l}{$0.330$} & \multicolumn{1}{l}{$0.206$}\\\hline
Analytical &  & \multicolumn{1}{l}{$0.1200$} & \multicolumn{1}{l}{$0.00288$} &
\multicolumn{1}{l}{$0.00006$} & \multicolumn{1}{l}{} & \multicolumn{1}{l}{}%
\\\hline
\end{tabular}
\ \ \ \ \ \ \ \ \ \ }\end{table}\begin{table}[ptb]
\caption{Comparison of $\ $the calculated invariants and error norms of
Problem 1 for $h=0.03$ and $k=0.2$ ($\mu=1,$ $3c=0.3,$ $x_{0}=10,$ $0\leq
x\leq30,$ $0\leq t\leq40$).}%
\label{p1t3}%
{\scriptsize
\begin{tabular}
[c]{ccccccc}\hline
Method & $t$ & $I_{1}$ & $I_{2}$ & $I_{3}$ & $L_{2}\times10^{3}$ & $L_{\infty
}\times10^{3}$\\\hline
CHCM-L & $0$ & \multicolumn{1}{l}{$1.1999445724$} &
\multicolumn{1}{l}{$0.2880000252$} & \multicolumn{1}{l}{$0.0576000000$} &
\multicolumn{1}{l}{$0.000679$} & \multicolumn{1}{l}{$0.003911$}\\
& $5$ & \multicolumn{1}{l}{$1.1999874409$} & \multicolumn{1}{l}{$0.2880000274$%
} & \multicolumn{1}{l}{$0.0576000012$} & \multicolumn{1}{l}{$0.000679$} &
\multicolumn{1}{l}{$0.019904$}\\
& $10$ & \multicolumn{1}{l}{$1.2000134431$} &
\multicolumn{1}{l}{$0.2880000287$} & \multicolumn{1}{l}{$0.0576000016$} &
\multicolumn{1}{l}{$0.025071$} & \multicolumn{1}{l}{$0.032148$}\\
& $20$ & \multicolumn{1}{l}{$1.2000387673$} &
\multicolumn{1}{l}{$0.2880000301$} & \multicolumn{1}{l}{$0.0576000018$} &
\multicolumn{1}{l}{$0.036187$} & \multicolumn{1}{l}{$0.044079$}\\
& $40$ & \multicolumn{1}{l}{$1.2000512978$} &
\multicolumn{1}{l}{$0.2880000311$} & \multicolumn{1}{l}{$0.0576000019$} &
\multicolumn{1}{l}{$0.048631$} & \multicolumn{1}{l}{$0.050084$}\\
CHCM-C & $40$ & $1.2000513180$ & $0.2880000297$ &
\multicolumn{1}{l}{$0.0576000015$} & \multicolumn{1}{l}{$0.052210$} &
\multicolumn{1}{l}{$0.050194$}\\
\cite{bl3zaki} & $40$ & \multicolumn{1}{l}{$1.1967$} &
\multicolumn{1}{l}{$0.2860$} & \multicolumn{1}{l}{$0.0570$} &
\multicolumn{1}{l}{$3.475$} & \multicolumn{1}{l}{$2.136$}\\
\cite{bl3kaslan} & $40$ & \multicolumn{1}{l}{$1.199992$} &
\multicolumn{1}{l}{$0.2921585$} & \multicolumn{1}{l}{$0.05759999$} &
\multicolumn{1}{l}{$0.07954512$} & \multicolumn{1}{l}{$-$}\\\hline
Analytical &  & \multicolumn{1}{l}{$1.2$} & \multicolumn{1}{l}{$0.288$} &
\multicolumn{1}{l}{$0.0576$} & \multicolumn{1}{l}{} & \multicolumn{1}{l}{}%
\\\hline
\end{tabular}
\ \ \ \ \ \ \ \ \ \ }\end{table}\begin{table}[ptb]
\caption{Comparison of $\ $the calculated invariants and error norms of
Problem 1 for various values of $N$ and $k$ at $t=40$ ($\mu=1,$ $3c=0.9,$
$x_{0}=40,$ $0\leq x\leq100$). }%
\label{p1t4}
{\scriptsize
\begin{tabular}
[c]{lcccccc}\hline
Method & $(N,k)$ & $I_{1}$ & $I_{2}$ & $I_{3}$ & $L_{2}$ &
\multicolumn{1}{c|}{$L_{\infty}$}\\\hline
CHCM-L &  & \multicolumn{1}{l}{} & \multicolumn{1}{l}{} & \multicolumn{1}{l}{}
& \multicolumn{1}{l}{} & \multicolumn{1}{l}{}\\
& \multicolumn{1}{l}{$(400,0.2)$} & \multicolumn{1}{l}{$3.5999999590$} &
\multicolumn{1}{l}{$2.5920297204$} & \multicolumn{1}{l}{$1.5552059235$} &
\multicolumn{1}{l}{$0.002671$} & \multicolumn{1}{l}{$0.001425$}\\
& \multicolumn{1}{l}{$(400,0.1)$} & \multicolumn{1}{l}{$3.5999999590$} &
\multicolumn{1}{l}{$2.5920296695$} & \multicolumn{1}{l}{$1.5552058615$} &
\multicolumn{1}{l}{$0.000696$} & \multicolumn{1}{l}{$0.000370$}\\
& \multicolumn{1}{l}{$(400,0.05)$} & \multicolumn{1}{l}{$3.5999999590$} &
\multicolumn{1}{l}{$2.5920296718$} & \multicolumn{1}{l}{$1.5552058581$} &
\multicolumn{1}{l}{$0.000202$} & \multicolumn{1}{l}{$0.000107$}\\
& \multicolumn{1}{l}{$(400,0.025)$} & \multicolumn{1}{l}{$3.5999999589$} &
\multicolumn{1}{l}{$2.5920296733$} & \multicolumn{1}{l}{$1.5552058580$} &
\multicolumn{1}{l}{$0.000079$} & \multicolumn{1}{l}{$0.000043$}\\
& \multicolumn{1}{l}{$(200,0.1)$} & \multicolumn{1}{l}{$3.5999974428$} &
\multicolumn{1}{l}{$2.5924441834$} & \multicolumn{1}{l}{$1.5552857179$} &
\multicolumn{1}{l}{$0.001251$} & \multicolumn{1}{l}{$0.000675$}\\
& \multicolumn{1}{l}{$(800,0.1)$} & \multicolumn{1}{l}{$3.5999999994$} &
\multicolumn{1}{l}{$2.5920018908$} & \multicolumn{1}{l}{$1.5552003783$} &
\multicolumn{1}{l}{$0.000661$} & \multicolumn{1}{l}{$0.000354$}\\
& \multicolumn{1}{l}{$(1600,0.1)$} & \multicolumn{1}{l}{$3.6000000000$} &
\multicolumn{1}{l}{$2.5920001237$} & \multicolumn{1}{l}{$1.5552000278$} &
\multicolumn{1}{l}{$0.000659$} & \multicolumn{1}{l}{$0.000353$}\\
& \multicolumn{1}{l}{} & \multicolumn{1}{l}{} & \multicolumn{1}{l}{} &
\multicolumn{1}{l}{} & \multicolumn{1}{l}{} & \multicolumn{1}{l}{}\\
CHCM-C & \multicolumn{1}{l}{} & \multicolumn{1}{l}{} & \multicolumn{1}{l}{} &
\multicolumn{1}{l}{} & \multicolumn{1}{l}{} & \multicolumn{1}{l}{}\\
& \multicolumn{1}{l}{$(400,0.2)$} & \multicolumn{1}{l}{$3.5998868279$} &
\multicolumn{1}{l}{$2.5918782898$} & \multicolumn{1}{l}{$1.5550774465$} &
\multicolumn{1}{l}{$0.004919$} & \multicolumn{1}{l}{$0.002894$}\\
& \multicolumn{1}{l}{$(400,0.1)$} & \multicolumn{1}{l}{$3.5998866685$} &
\multicolumn{1}{l}{$2.5918790458$} & \multicolumn{1}{l}{$1.5550779337$} &
\multicolumn{1}{l}{$0.003078$} & \multicolumn{1}{l}{$0.001896$}\\
& \multicolumn{1}{l}{$(400,0.05)$} & \multicolumn{1}{l}{$3.5998866286$} &
\multicolumn{1}{l}{$2.5918792500$} & \multicolumn{1}{l}{$1.5550780676$} &
\multicolumn{1}{l}{$0.002647$} & \multicolumn{1}{l}{$0.001646$}\\
& \multicolumn{1}{l}{$(400,0.025)$} & \multicolumn{1}{l}{$3.5998866187$} &
\multicolumn{1}{l}{$2.5918793020$} & \multicolumn{1}{l}{$1.5550781018$} &
\multicolumn{1}{l}{$0.002542$} & \multicolumn{1}{l}{$0.001584$}\\
& \multicolumn{1}{l}{$(200,0.1)$} & \multicolumn{1}{l}{$3.5982884438$} &
\multicolumn{1}{l}{$2.5901738078$} & \multicolumn{1}{l}{$1.553355029$} &
\multicolumn{1}{l}{$0.011636$} & \multicolumn{1}{l}{$0.007023$}\\
& \multicolumn{1}{l}{$(800,0.1)$} & \multicolumn{1}{l}{$3.5999928133$} &
\multicolumn{1}{l}{$2.5919922760$} & \multicolumn{1}{l}{$1.5551922231$} &
\multicolumn{1}{l}{$0.001212$} & \multicolumn{1}{l}{$0.000717$}\\
& \multicolumn{1}{l}{$(1600,0.1)$} & \multicolumn{1}{l}{$3.5999995492$} &
\multicolumn{1}{l}{$2.5919995049$} & \multicolumn{1}{l}{$1.5551995049$} &
\multicolumn{1}{l}{$0.000788$} & \multicolumn{1}{l}{$0.000433$}\\
& \multicolumn{1}{l}{} & \multicolumn{1}{l}{} & \multicolumn{1}{l}{} &
\multicolumn{1}{l}{} & \multicolumn{1}{l}{} & \multicolumn{1}{l}{}\\
\cite{ek2} & \multicolumn{1}{l}{} & \multicolumn{1}{l}{} &
\multicolumn{1}{l}{} & \multicolumn{1}{l}{} & \multicolumn{1}{l}{} &
\multicolumn{1}{l}{}\\
& \multicolumn{1}{l}{$(400,0.2)$ EXE} & \multicolumn{1}{l}{$3.600000$} &
\multicolumn{1}{l}{$2.882298$} & \multicolumn{1}{l}{$1.828214$} &
\multicolumn{1}{l}{$0.0133293$} & $-$\\
& \multicolumn{1}{l}{$(400,0.1)$ EXE} & \multicolumn{1}{l}{$3.599999$} &
\multicolumn{1}{l}{$2.724104$} & \multicolumn{1}{l}{$1.681432$} &
\multicolumn{1}{l}{$0.00490421$} & $-$\\
& \multicolumn{1}{l}{$(400,0.05)$ EXE} & \multicolumn{1}{l}{$3.600000$} &
\multicolumn{1}{l}{$2.652641$} & \multicolumn{1}{l}{$1.616028$} &
\multicolumn{1}{l}{$0.00247959$} & $-$\\
& \multicolumn{1}{l}{$(400,0.025)$ E} & \multicolumn{1}{l}{$3.600000$} &
\multicolumn{1}{l}{$2.652160$} & \multicolumn{1}{l}{$1.615533$} &
\multicolumn{1}{l}{$0.00310718$} & $-$\\
& \multicolumn{1}{l}{$(200,0.1)$ EXE} & \multicolumn{1}{l}{$3.600000$} &
\multicolumn{1}{l}{$2.837960$} & \multicolumn{1}{l}{$1.807345$} &
\multicolumn{1}{l}{$0.0105221$} & $-$\\
& \multicolumn{1}{l}{$(800,0.1)$ E} & \multicolumn{1}{l}{$3.600000$} &
\multicolumn{1}{l}{$2.893544$} & \multicolumn{1}{l}{$1.833380$} &
\multicolumn{1}{l}{$0.0163510$} & $-$\\
& \multicolumn{1}{l}{$(1600,0.1)$ E} & \multicolumn{1}{l}{$3.600000$} &
\multicolumn{1}{l}{$2.896941$} & \multicolumn{1}{l}{$1.835213$} &
\multicolumn{1}{l}{$0.0173992$} & $-$\\\hline
\end{tabular}
\ \ \ \ \ \ \ \ \ \ }\end{table}\begin{table}[ptb]
\caption{Comparison of $\ $the calculated invariants and error norms of
Problem 1 for various values of $N$ and $k$ at $t=40$ ($\mu=1,$ $3c=0.9,$
$x_{0}=40,$ $0\leq x\leq100$).
\ \ \ \ \ \ \ \ \ \ \ \ \ \ \ \ \ \ \ \ \ \ \ \ \ \ \ \ \ \ \ \ \ \ \ \ \ \ \ \ \ \ \ \ \ }%
\label{p1t5}
{\scriptsize
\begin{tabular}
[c]{ccccccc}\hline
Method & $(N,k)$ & $I_{1}$ & $I_{2}$ & $I_{3}$ & $L_{2}\times10^{3}$ &
\multicolumn{1}{c|}{$L_{\infty}\times10^{3}$}\\\hline
CHCM-L & \multicolumn{1}{l}{} & \multicolumn{1}{l}{} & \multicolumn{1}{l}{} &
\multicolumn{1}{l}{} & \multicolumn{1}{l}{} & \multicolumn{1}{l}{}\\
& \multicolumn{1}{l}{$(400,0.01)$} & \multicolumn{1}{l}{$3.5999999589$} &
\multicolumn{1}{l}{$2.5920296738$} & \multicolumn{1}{l}{$1.5552058580$} &
\multicolumn{1}{l}{$0.045036$} & \multicolumn{1}{l}{$0.025156$}\\
& \multicolumn{1}{l}{$(400,0.005)$} & \multicolumn{1}{l}{$3.5999999590$} &
\multicolumn{1}{l}{$2.5920296739$} & \multicolumn{1}{l}{$1.5552058580$} &
\multicolumn{1}{l}{$0.040248$} & \multicolumn{1}{l}{$0.022653$}\\
& \multicolumn{1}{l}{$(400,0.0025)$} & \multicolumn{1}{l}{$3.5999999590$} &
\multicolumn{1}{l}{$2.5920296739$} & \multicolumn{1}{l}{$1.5552058580$} &
\multicolumn{1}{l}{$0.039056$} & \multicolumn{1}{l}{$0.022027$}\\
& \multicolumn{1}{l}{$(400,0.00125)$} & \multicolumn{1}{l}{$3.5999999590$} &
\multicolumn{1}{l}{$2.5920296739$} & \multicolumn{1}{l}{$1.5552058580$} &
\multicolumn{1}{l}{$0.038759$} & \multicolumn{1}{l}{$0.021871$}\\
& \multicolumn{1}{l}{$(800,0.000625)$} & \multicolumn{1}{l}{$3.5999999994$} &
\multicolumn{1}{l}{$2.5920018861$} & \multicolumn{1}{l}{$1.5552003741$} &
\multicolumn{1}{l}{$0.002455$} & \multicolumn{1}{l}{$0.001389$}\\
& \multicolumn{1}{l}{$(1600,0.0003125)$} & \multicolumn{1}{l}{$3.6000000000$}
& \multicolumn{1}{l}{$2.5920001184$} & \multicolumn{1}{l}{$1.5552000235$} &
\multicolumn{1}{l}{$0.000158$} & \multicolumn{1}{l}{$0.000089$}\\
&  & \multicolumn{1}{l}{} & \multicolumn{1}{l}{} & \multicolumn{1}{l}{} &
\multicolumn{1}{l}{} & \multicolumn{1}{l}{}\\
CHCM-C & \multicolumn{1}{l}{} & \multicolumn{1}{l}{} & \multicolumn{1}{l}{} &
\multicolumn{1}{l}{} & \multicolumn{1}{l}{} & \multicolumn{1}{l}{}\\
& \multicolumn{1}{l}{$(400,0.01)$} & \multicolumn{1}{l}{$3.5998866159$} &
\multicolumn{1}{l}{$2.5918793166$} & \multicolumn{1}{l}{$1.5550781115$} &
\multicolumn{1}{l}{$2.513429$} & \multicolumn{1}{l}{$1.566399$}\\
& \multicolumn{1}{l}{$(400,0.005)$} & \multicolumn{1}{l}{$3.5998866155$} &
\multicolumn{1}{l}{$2.5918793187$} & \multicolumn{1}{l}{$1.5550781128$} &
\multicolumn{1}{l}{$2.509307$} & \multicolumn{1}{l}{$1.563904$}\\
& \multicolumn{1}{l}{$(400,0.0025)$} & \multicolumn{1}{l}{$3.5998866154$} &
\multicolumn{1}{l}{$2.5918793193$} & \multicolumn{1}{l}{$1.5550781132$} &
\multicolumn{1}{l}{$2.507277$} & \multicolumn{1}{l}{$1.563280$}\\
& \multicolumn{1}{l}{$(400,0.00125)$} & \multicolumn{1}{l}{$3.5998866154$} &
\multicolumn{1}{l}{$2.5918793194$} & \multicolumn{1}{l}{$1.5550781133$} &
\multicolumn{1}{l}{$2.508020$} & \multicolumn{1}{l}{$1.563124$}\\
& \multicolumn{1}{l}{$(800,0.000625)$} & \multicolumn{1}{l}{$3.5999928099$} &
\multicolumn{1}{l}{$2.5919923510$} & \multicolumn{1}{l}{$1.5551922755$} &
\multicolumn{1}{l}{$0.610566$} & \multicolumn{1}{l}{$0.384542$}\\
& \multicolumn{1}{l}{$(1600,0.0003125)$} & \multicolumn{1}{l}{$3.5999995490$}
& \multicolumn{1}{l}{$2.5919995203$} & \multicolumn{1}{l}{$1.5551995156$} &
\multicolumn{1}{l}{$0.151612$} & \multicolumn{1}{l}{$0.095748$}\\
&  & \multicolumn{1}{l}{} & \multicolumn{1}{l}{} & \multicolumn{1}{l}{} &
\multicolumn{1}{l}{} & \multicolumn{1}{l}{}\\
\cite{ek2}CE & \multicolumn{1}{l}{} & \multicolumn{1}{l}{} &
\multicolumn{1}{l}{} & \multicolumn{1}{l}{} & \multicolumn{1}{l}{} &
\multicolumn{1}{l}{}\\
& \multicolumn{1}{l}{$(400,0.01)$} & \multicolumn{1}{l}{$3.599999$} &
\multicolumn{1}{l}{$2.612544$} & \multicolumn{1}{l}{$1.579549$} &
\multicolumn{1}{l}{$1.61169$} & $-$\\
& \multicolumn{1}{l}{$(400,0.005)$} & \multicolumn{1}{l}{$3.600000$} &
\multicolumn{1}{l}{$2.599007$} & \multicolumn{1}{l}{$1.567285$} &
\multicolumn{1}{l}{$0.795278$} & $-$\\
& \multicolumn{1}{l}{$(400,0.0025)$} & \multicolumn{1}{l}{$3.600000$} &
\multicolumn{1}{l}{$2.592304$} & \multicolumn{1}{l}{$1.561220$} &
\multicolumn{1}{l}{$0.389948$} & $-$\\
& \multicolumn{1}{l}{$(400,0.00125)$} & \multicolumn{1}{l}{$3.600001$} &
\multicolumn{1}{l}{$2.588970$} & \multicolumn{1}{l}{$1.558205$} &
\multicolumn{1}{l}{$0.188069$} & $-$\\
& \multicolumn{1}{l}{$(800,0.000625)$} & \multicolumn{1}{l}{$3.599999$} &
\multicolumn{1}{l}{$2.592064$} & \multicolumn{1}{l}{$1.556701$} &
\multicolumn{1}{l}{$0.0999448$} & $-$\\
& \multicolumn{1}{l}{$(1600,0.0003125)$} & \multicolumn{1}{l}{$3.600000$} &
\multicolumn{1}{l}{$2.592432$} & \multicolumn{1}{l}{$1.555950$} &
\multicolumn{1}{l}{$0.0503415$} & $-$\\\hline
\end{tabular}
\ \ \ \ \ \ \ \ \ \ }\end{table}The graphs of the simulations of single
solitary wave for different values of velocity and amplitudes are plotted in
Figure $\ref{F1}.$ One can easily see from Figure $\ref{F1}$ that the
amplitudes, velocities and shapes of the wave are conserved during the
simulation. Furthermore, in Table \ref{p1t1}, one can see the comparison of
our results with some of those existing in the literature. From the table, it
is observed that\ the newly obtained results\ are better than the other ones
except those given in Refs. \cite{bl3bsaka} and \cite{bl3esen}. Table
$\ref{p1t2}$ shows a comparison of $\ $the $3$ invariants and the error norms
of Problem 1 for $h=k=0.05$ ($\mu=1,$ $3c=0.03,$ $x_{0}=10,$ $0\leq x\leq30,$
$0\leq t\leq80$) with their analytical values and those in Refs.
\cite{bl3esen} and \cite{bl3dogan}. Again Table $\ref{p1t3}$ presents a clear
comparison of $\ $the $3$ invariants and the error norms of Problem 1 for
values of $k=0.2$ and $h=0.03$ ($3c=0.3,$ $\mu=1,$ $x_{0}=10,$ $0\leq
x\leq30,$ $0\leq t\leq40$) with their analytical values and those in Refs.
\cite{bl3zaki} and \cite{bl3kaslan}.

Table \ref{p1t4} shows a comparison of $\ $the $3$ invariants and also the
error norms of Problem 1 for various values of $N$ and $k$ at $t=40$ ($\mu=1,$
$3c=0.9,$ $x_{0}=40,$ $0\leq x\leq100$). One can clearly see from Table
\ref{p1t4} that those results found by taking the shifted roots of the
Legendre polynomial as interior collocation points required in the proposed
method are much better than those obtained by taking the shifted roots of the
Chebyshev polynomial. Since it is known that Legendre polynomials minimize the
error and give appropriate results, such results were expected beforehand.
Finally, Table \ref{p1t5} shows a comparison of $\ $the $3$ invariants and
also the error norms of Problem 1 for various values of $N$ and $k$ at $t=40$
($\mu=1,$ $3c=0.9,$ $x_{0}=40,$ $0\leq x\leq100$). One can also obviously see
from both Tables \ref{p1t4} and \ref{p1t5} as $h$ and $k$ decrease so the
values of the error norms $L_{2}$ and $L_{\infty}$ decrease. In other words,
the obtained numerical solution approaches to the analytical solution. This
shows that numerical solutions satisfy the expected
accuracy.\begin{figure}[ptb]
\begin{center}
\centering\includegraphics[width=0.48\textwidth ]{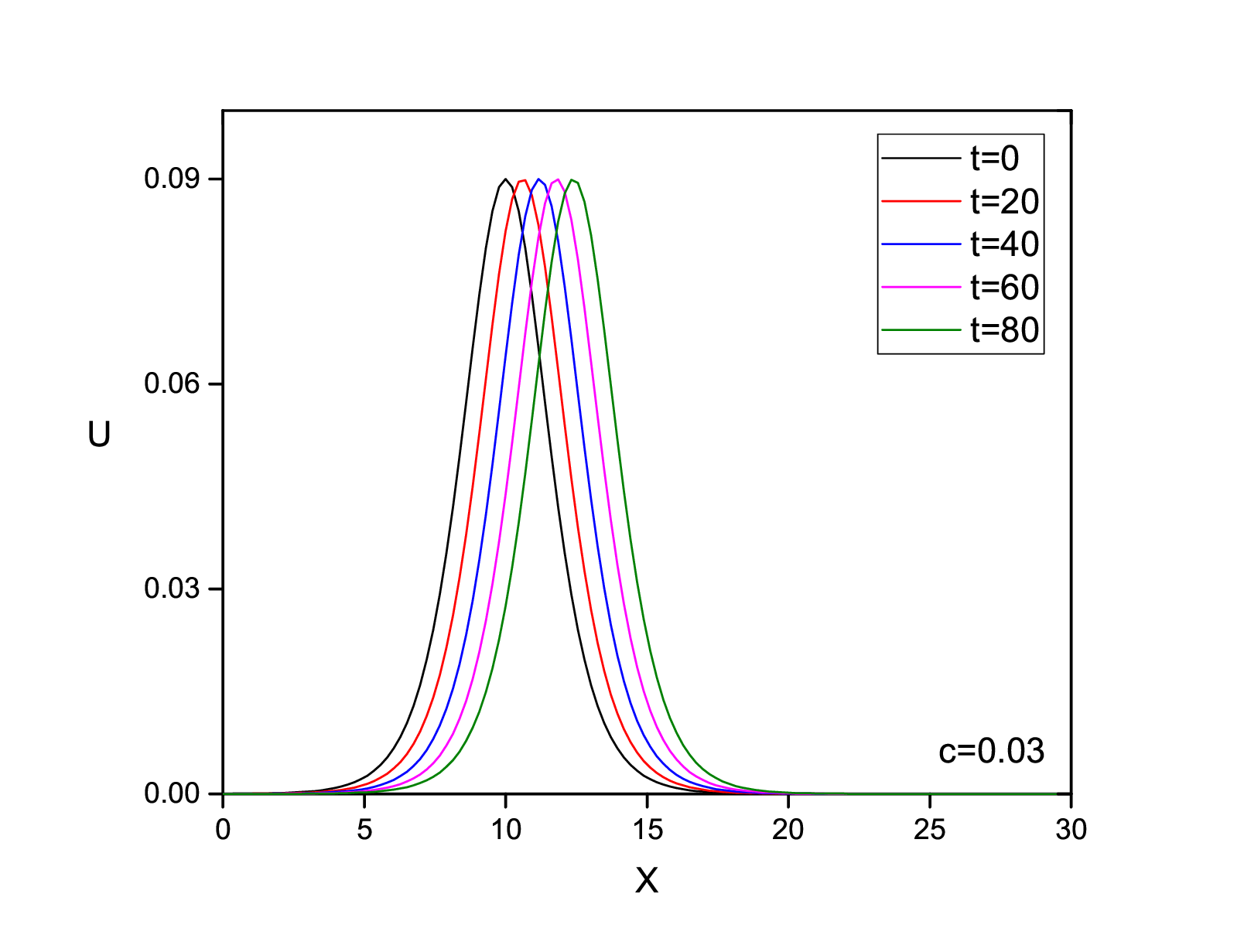}
\includegraphics[width=0.48\textwidth ]{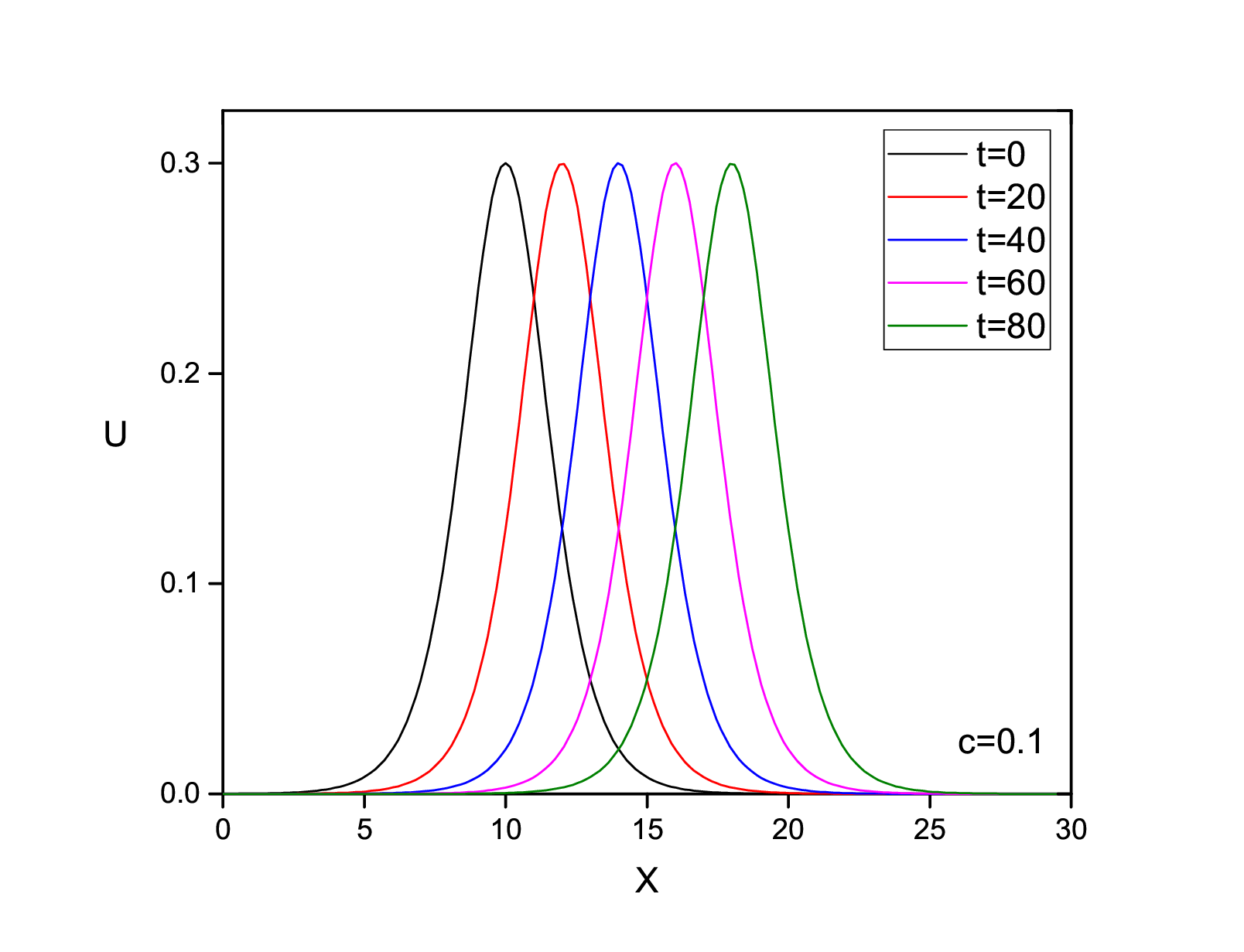}
\end{center}
\caption{\textrm{{\protect\footnotesize Simulations} {\protect\footnotesize of
single solitary wave for velocity values} }{\protect\footnotesize c=}%
${\protect\footnotesize 0.03}${\protect\footnotesize , }%
${\protect\footnotesize 0.1}$\textrm{\ }{\protect\footnotesize at
}${\protect\footnotesize t=0(20)80}${\protect\footnotesize .}}%
\label{F1}%
\end{figure}

\subsection{Two solitary waves}

The second test experimental problem has been taken as the interaction of $2$
solitary waves. We are going to take into consideration Eq. $(\ref{1})$ with
the solution domain $(x,t)\in\lbrack a,b]\times\lbrack0,T],$ the initial
condition \cite{bl3esen} and the boundary conditions (\ref{bc})%

\[
U(x,0)=\overset{2}{\underset{j=1}{\sum}}3c_{j}\sec h^{2}\left[  0.5\left(
x-x_{j}-c_{j}\right)  \right]
\]
where the parameters $\mu=1,c_{1}=1.5$, $c_{2}=0.75$, $x_{1}=10,$ $x_{2}=25$
with $\Delta t=0.01$ are taken in the region $0\leq x\leq80.$ The exact values
of the invariants are found as $I_{1}=12\left(  c_{1}+c_{2}\right)  =27,$
$I_{2}=28.8\left(  c_{1}^{2}+c_{2}^{2}\right)  =81\ $and $I_{3}=57.6\left(
c_{1}^{3}+c_{2}^{3}\right)  =218.7.$

\begin{figure}[ptb]
\begin{center}
\centering\includegraphics[width=0.48\textwidth]{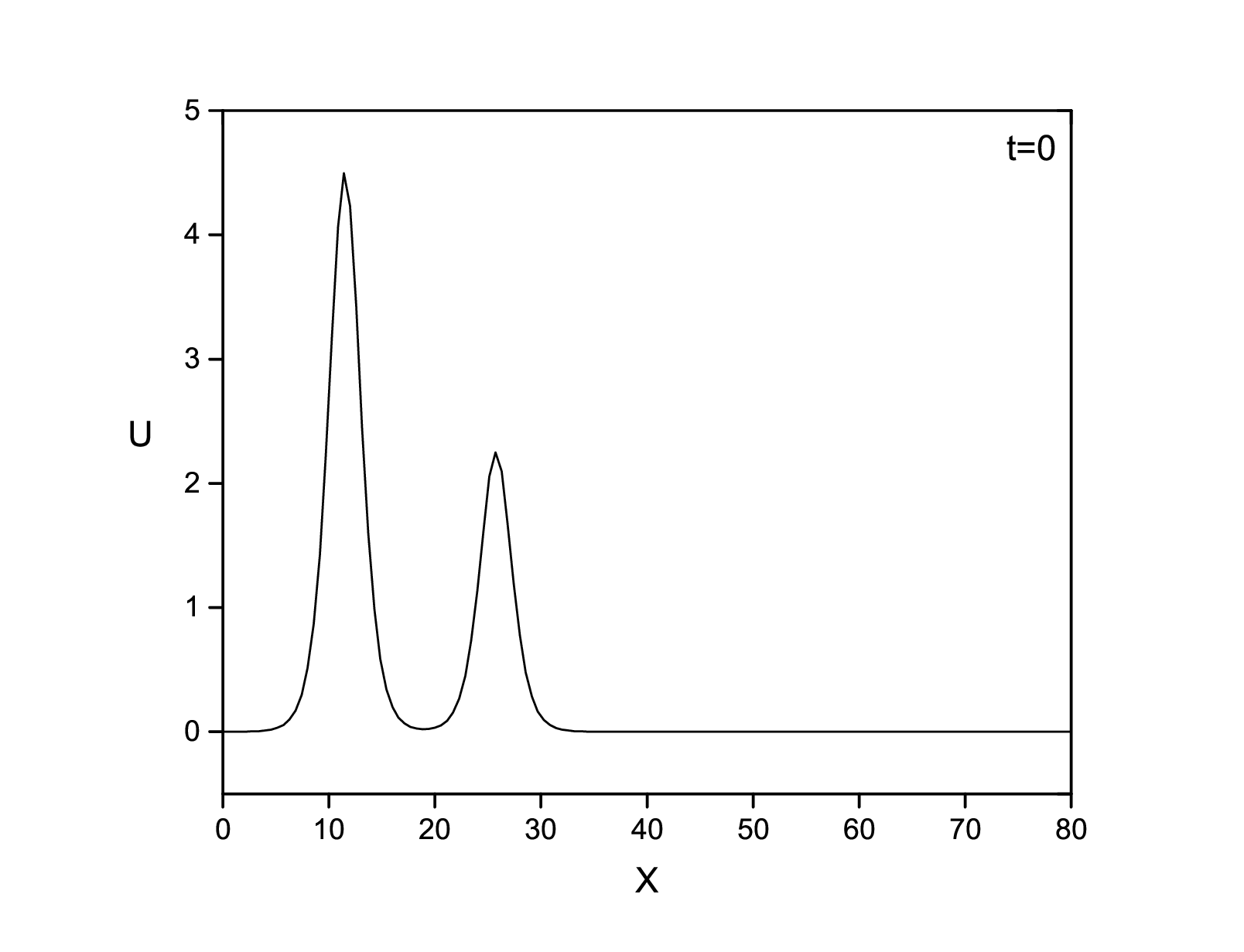}
\includegraphics[width=0.48\textwidth]{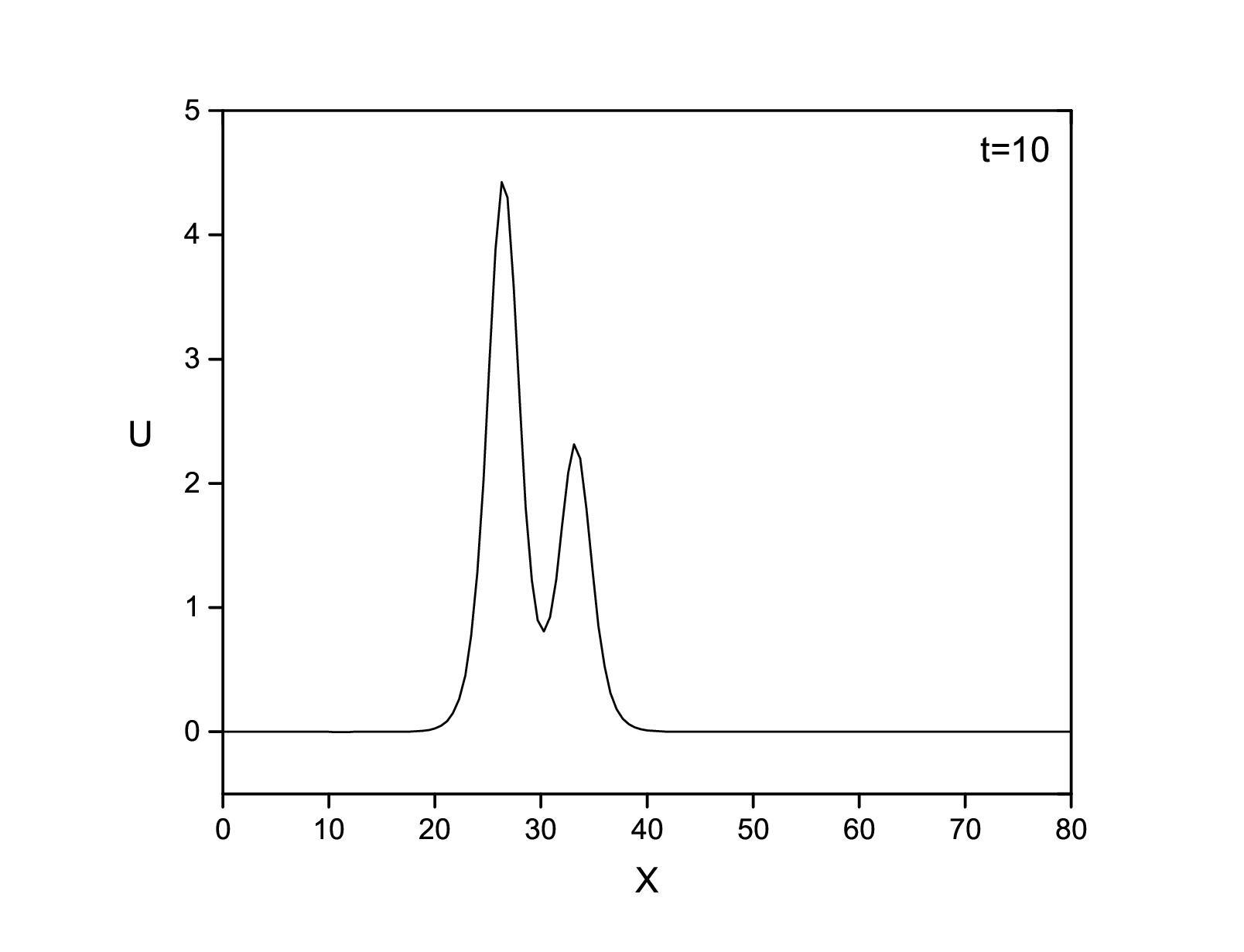}
\includegraphics[width=0.48\textwidth]{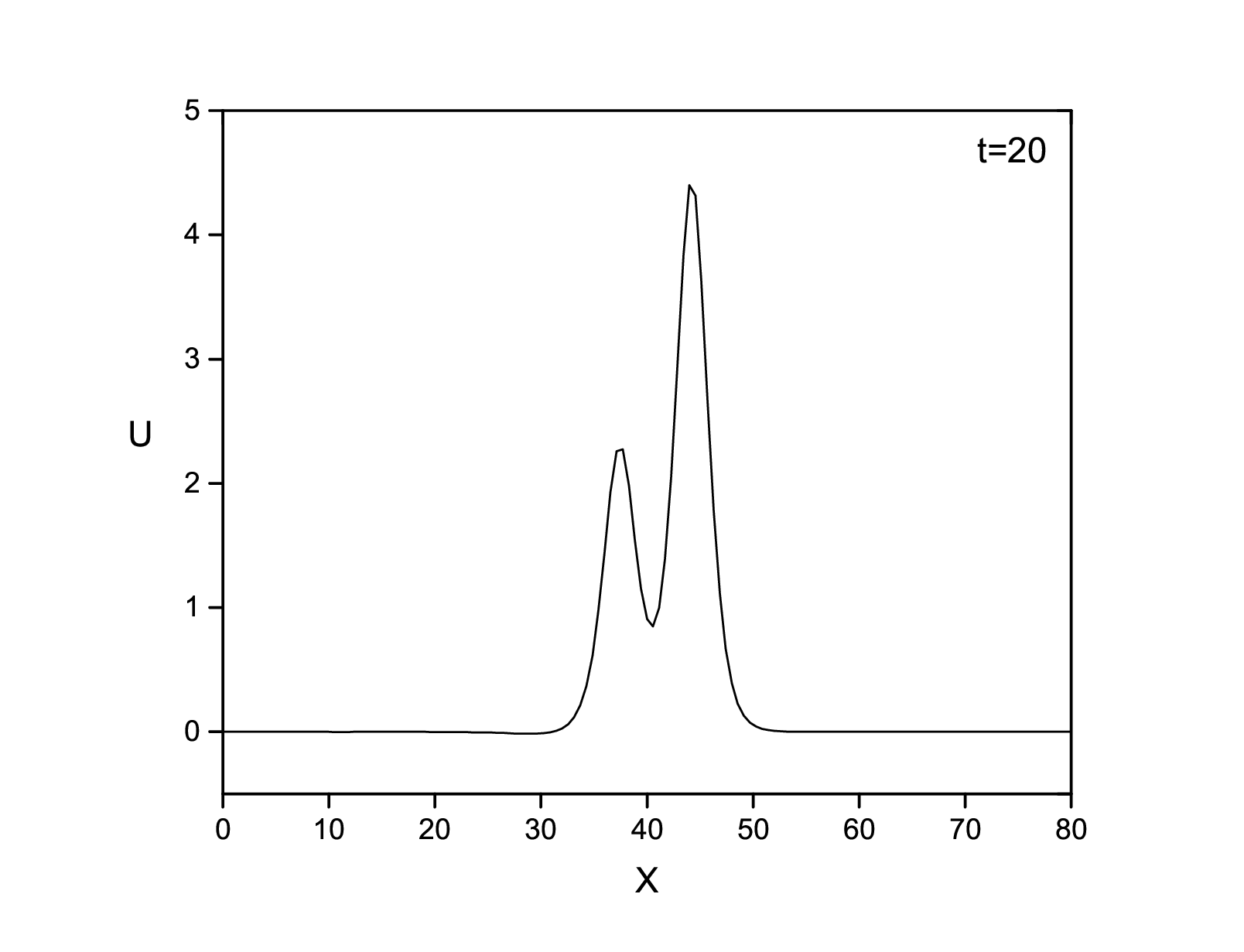}
\includegraphics[width=0.48\textwidth]{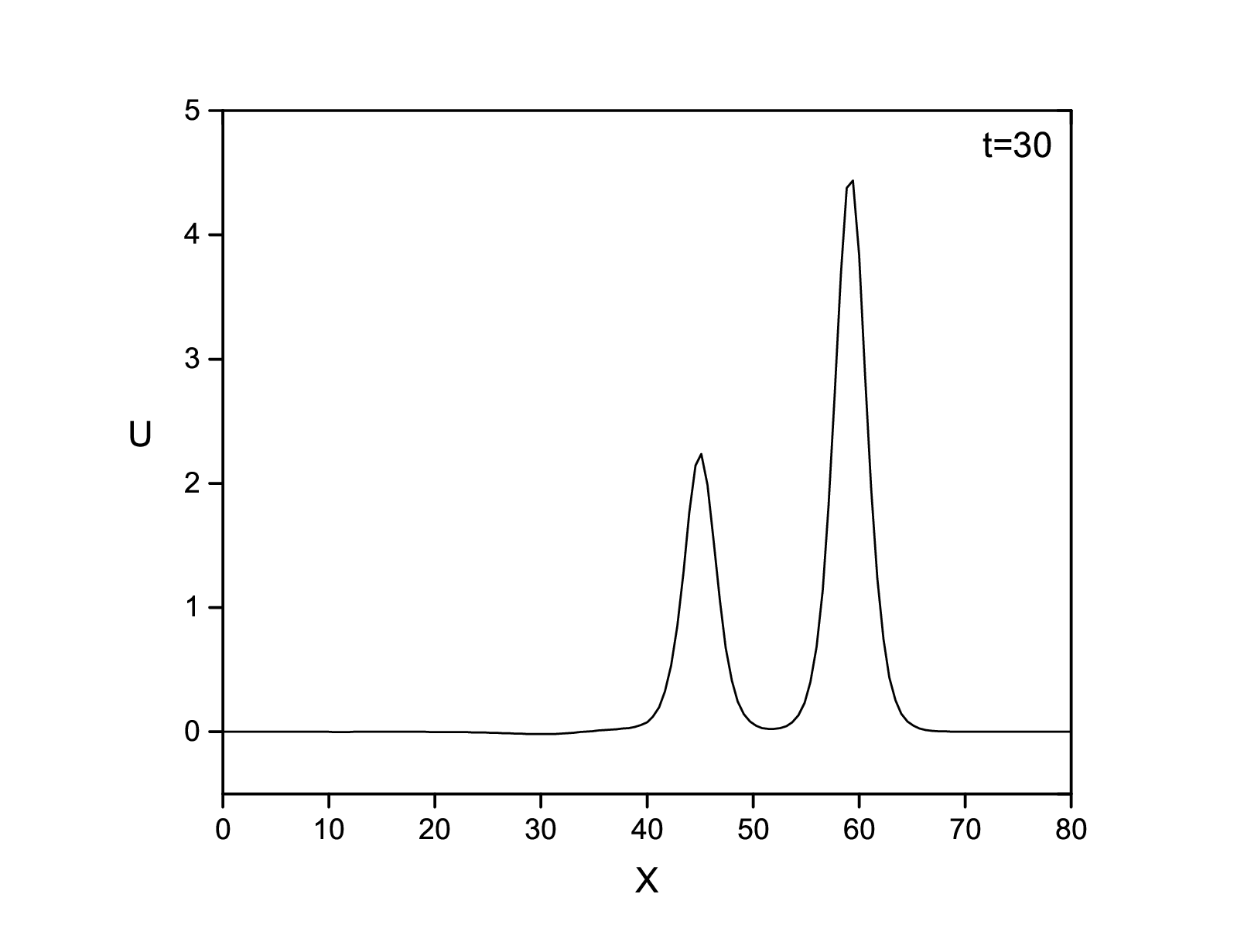}
\end{center}
\caption{The {\protect\footnotesize simula}%
\textrm{{\protect\footnotesize tion} {\protect\footnotesize of }%
}$\mathrm{{\protect\footnotesize 2}}$\textrm{{\protect\footnotesize solitary
waves at times }}$t=0,$ $10,$ $20,$ $30{\protect\footnotesize .}$}%
\label{F02}%
\end{figure}

The simulation of the interaction of $2$ solitary waves until time $t=30$ is
presented in Figure $\ref{F02}.$ One can easily see from this figure that the
interaction process started approximately at $t=10$ , and the separation
process started approximately at $t=20.$ In the end, $2$ waves replaced their
initial positions. In Table \ref{p2t1}, the calculated results have been
compared to those existing in the literature. One can obviously see from this
table that the newly obtained results are$\ $in good harmony with their exact
values and also all of the compared ones.\begin{table}[ptb]
\caption{Comparison of $\ $the calculated invariants of Problem 2 for
$h=k=0.1$ ($\mu=1,$ $c_{1}=1.5,$ $c_{2}=0.75,$ $x_{1}=10$, $x_{2}=25,$ $0\leq
x\leq80,0\leq t\leq30$ ).}%
\label{p2t1}%
{\scriptsize
\begin{tabular}
[c]{clcccc}\hline
Method &  & $t$ & $I_{1}$ & $I_{2}$ & $I_{3}$\\\hline
CHCM-L &  & $1$ & \multicolumn{1}{l}{$27.000090$} &
\multicolumn{1}{l}{$81.000450$} & \multicolumn{1}{l}{$218.702919$}\\
&  & $5$ & \multicolumn{1}{l}{$27.000171$} & \multicolumn{1}{l}{$81.000368$} &
\multicolumn{1}{l}{$218.702149$}\\
&  & $10$ & \multicolumn{1}{l}{$27.000171$} & \multicolumn{1}{l}{$80.994156$}
& \multicolumn{1}{l}{$218.662061$}\\
&  & $15$ & \multicolumn{1}{l}{$27.000171$} & \multicolumn{1}{l}{$80.940889$}
& \multicolumn{1}{l}{$218.323702$}\\
&  & $20$ & \multicolumn{1}{l}{$27.000171$} & \multicolumn{1}{l}{$80.992358$}
& \multicolumn{1}{l}{$218.653188$}\\
&  & $25$ & \multicolumn{1}{l}{$27.000171$} & \multicolumn{1}{l}{$81.000154$}
& \multicolumn{1}{l}{$218.701589$}\\
&  & $30$ & \multicolumn{1}{l}{$27.000171$} & \multicolumn{1}{l}{$81.000478$}
& \multicolumn{1}{l}{$218.703143$}\\
CHCM-C &  & $30$ & \multicolumn{1}{l}{$27.000066$} &
\multicolumn{1}{l}{$80.999907$} & \multicolumn{1}{l}{$218.700877$}\\
\cite{ek1} &  & $30$ & \multicolumn{1}{l}{$26.999994$} &
\multicolumn{1}{l}{$81.000511$} & \multicolumn{1}{l}{$218.703446$}\\
\cite{bl014} &  & $30$ & \multicolumn{1}{l}{$27.00017$} &
\multicolumn{1}{l}{$80.96848$} & \multicolumn{1}{l}{$218.70210$}\\
\cite{bl3esen} &  & $30$ & \multicolumn{1}{l}{$27.00003$} &
\multicolumn{1}{l}{$81.01719$} & \multicolumn{1}{l}{$218.70650$}\\
\cite{bl3bsaka} &  & $30$ & \multicolumn{1}{l}{$27.00068$} &
\multicolumn{1}{l}{$81.02407$} & \multicolumn{1}{l}{$218.73673$}\\
\cite{bl3fazal} &  & $30$ & \multicolumn{1}{l}{$27.00019$} &
\multicolumn{1}{l}{$81.00045$} & \multicolumn{1}{l}{$218.70312$}\\
\cite{bl3kraslan} &  & $30$ & \multicolumn{1}{l}{$27.12702$} &
\multicolumn{1}{l}{$80.98988$} & \multicolumn{1}{l}{$218.6996$}\\
\cite{bl3ali} & $(h=0.4)$ & $30$ & \multicolumn{1}{l}{$27.00000$} &
\multicolumn{1}{l}{$80.999703$} & \multicolumn{1}{l}{$218.69966$}\\
\cite{bl3bsakadag} &  & $30$ & \multicolumn{1}{l}{$27.00017$} &
\multicolumn{1}{l}{$81.00044$} & \multicolumn{1}{l}{$218.70304$}\\
\cite{bl3ydereli} & $(h=0.4)$ & $30$ & \multicolumn{1}{l}{$27.000582$} &
\multicolumn{1}{l}{$81.001095$} & \multicolumn{1}{l}{$218.726082$}\\
\cite{bl3uddin} & $(h=0.2,$ $k=0.05)$ & $30$ & \multicolumn{1}{l}{$26.93310$}
& \multicolumn{1}{l}{$80.80028$} & \multicolumn{1}{l}{$218.16659$}\\\hline
Analytical &  &  & \multicolumn{1}{l}{$27$} & \multicolumn{1}{l}{$81$} &
\multicolumn{1}{l}{$218.7$}\\\hline
\end{tabular}
\ \ \ \ \ \ \ \ \ \ }\end{table}

\subsection{Three solitary waves}

The third experimental problem is the interaction of $3$ solitary waves. Eq.
$(\ref{1})$ will be considered over solution domain $(x,t)\in\lbrack
a,b]\times\lbrack0,T],$ and the boundary conditions (\ref{bc}) and the initial
condition \cite{bl013}%

\[
U(x,0)=\overset{3}{\underset{j=1}{\sum}}3c_{j}\sec h^{2}\left[  0.5\left(
x-x_{j}-c_{j}\right)  \right]
\]
in which the parameters $\mu=1,c_{1}=4.5$, $c_{2}=1.5$, $c_{3}=0.5$,
$x_{1}=10,$ $x_{2}=25,$ $x_{3}=35$ with $\Delta t=0.1$ are taken over the
region $\left[  0,100\right]  .$ Therefore, the analytical values of the
invariants can be found as $I_{1}=12\left(  c_{1}+c_{2}+c_{3}\right)  =78,$
$I_{2}=28.8\left(  c_{1}^{2}+c_{2}^{2}+c_{3}^{2}\right)  =655.2$ and
$I_{3}=57.6\left(  c_{1}^{3}+c_{2}^{3}+c_{3}^{3}\right)  =5450.4.$

\begin{figure}[ptb]
\begin{center}
\centering\includegraphics[width=0.48\textwidth]{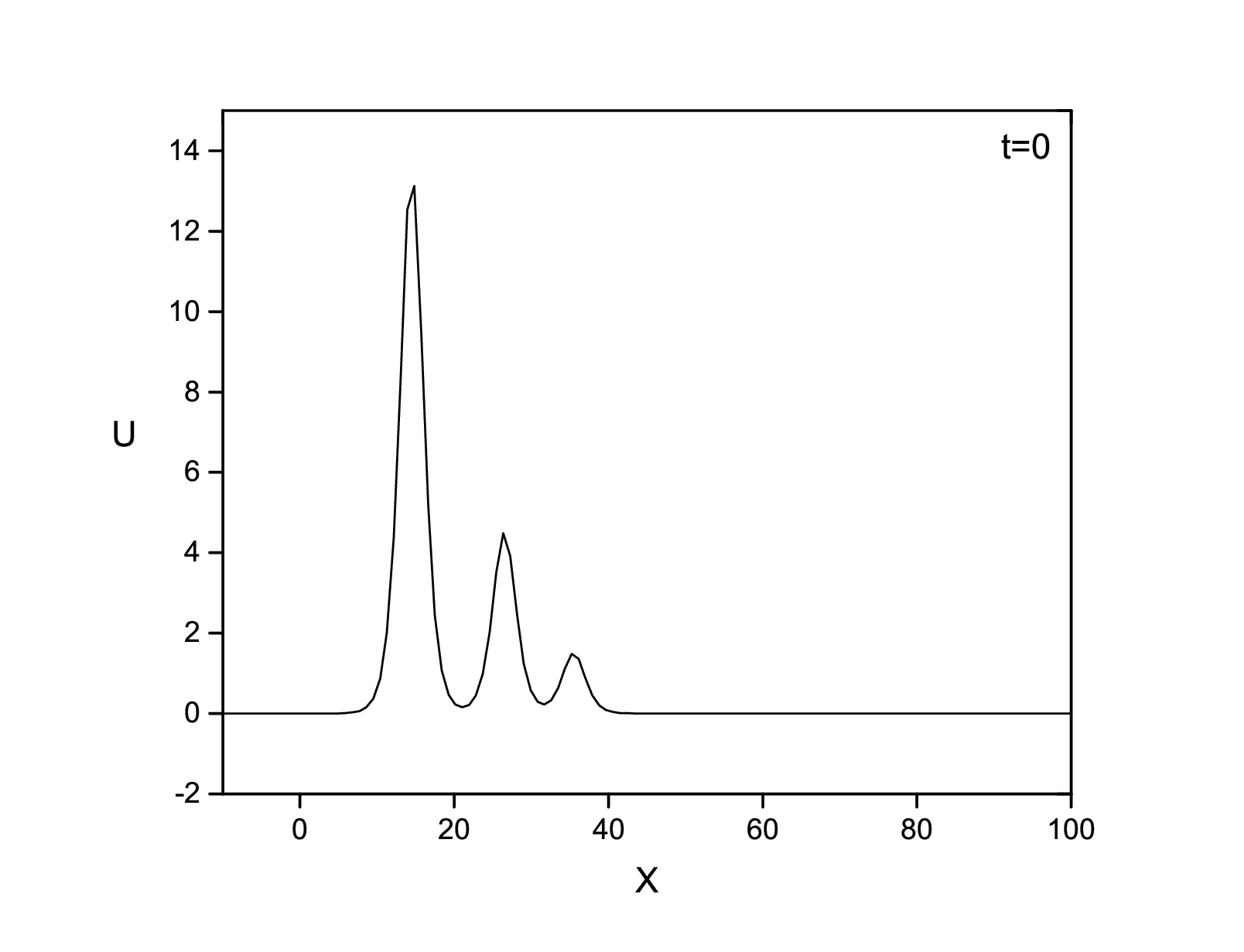}
\includegraphics[width=0.48\textwidth]{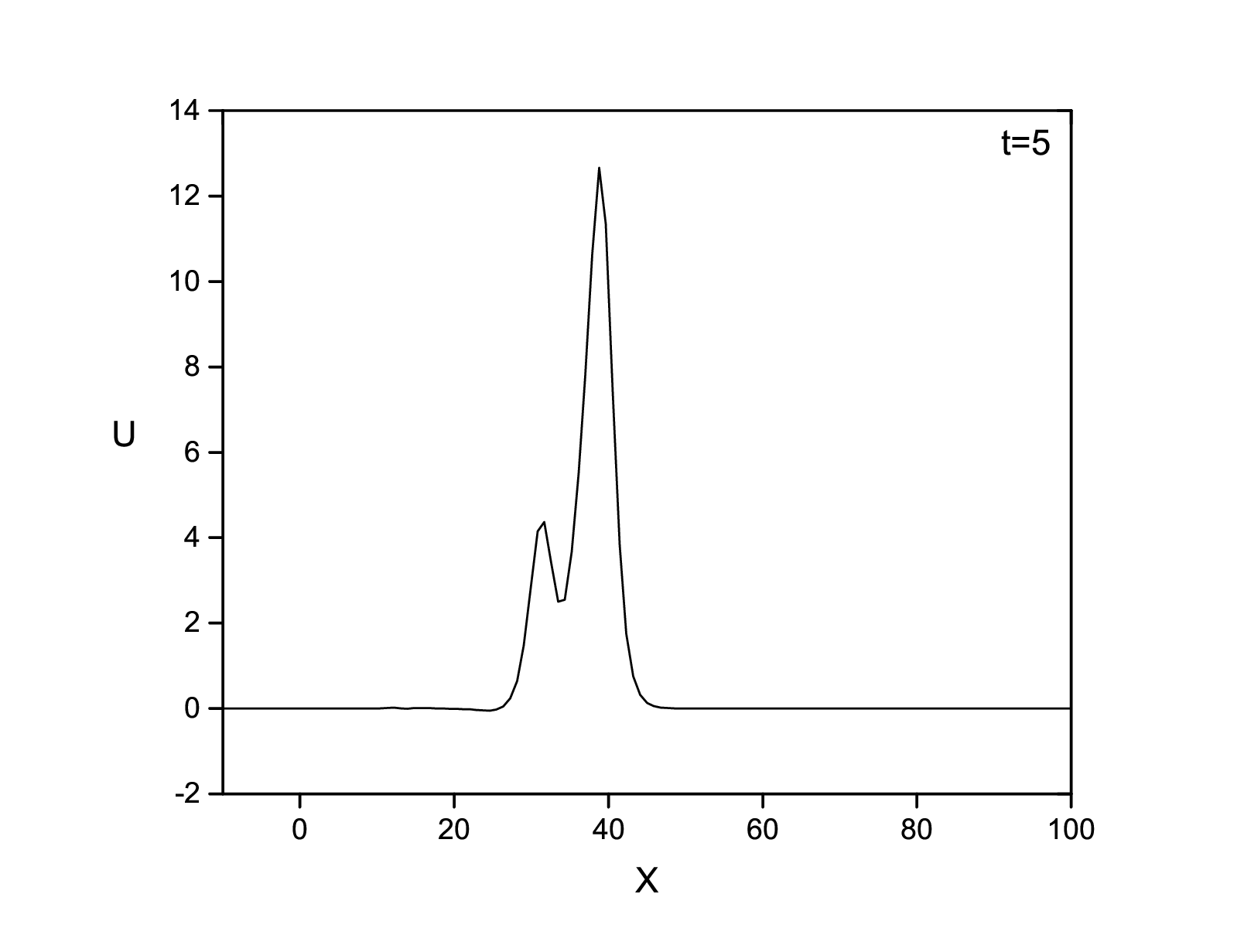}
\includegraphics[width=0.48\textwidth]{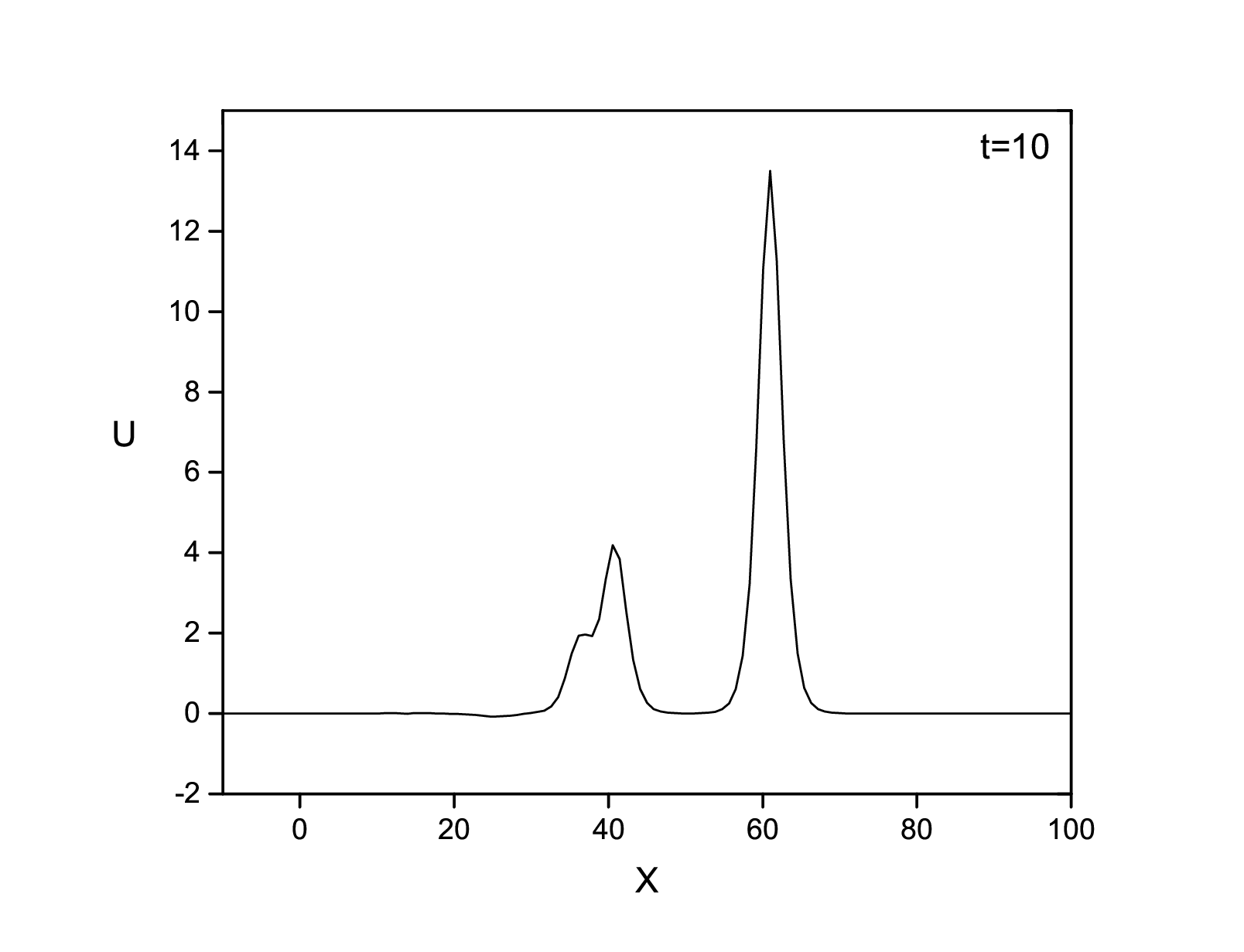}
\includegraphics[width=0.48\textwidth]{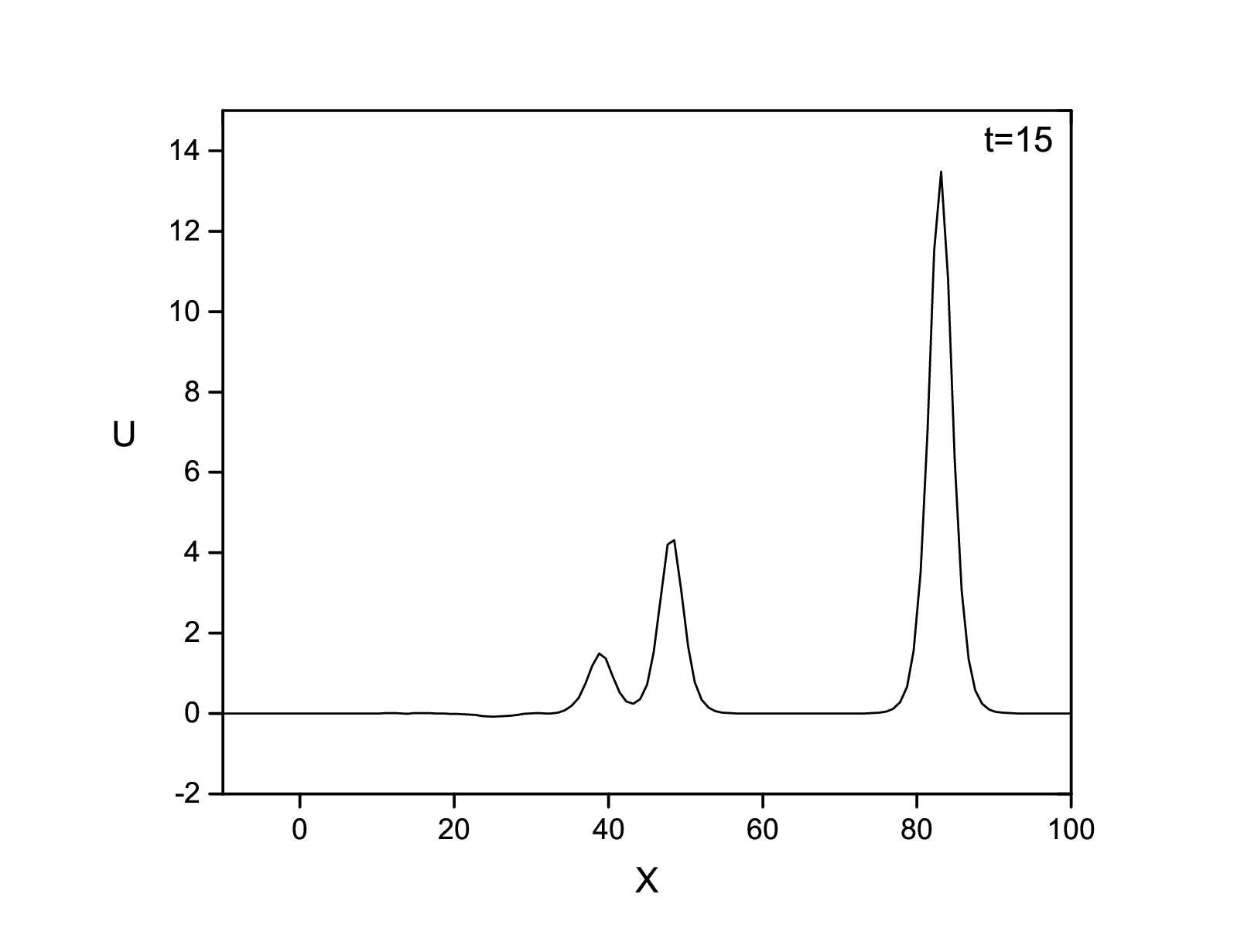}
\end{center}
\caption{The {\protect\footnotesize simula}%
\textrm{{\protect\footnotesize tion} {\protect\footnotesize of }%
}$\mathrm{{\protect\footnotesize 3}}$\textrm{{\protect\footnotesize solitary
waves at times }}$t=0,$ $5,$ $10,$ $15{\protect\footnotesize .}$}%
\label{F3}%
\end{figure}

The simulation for the interaction of $3$ solitary waves run up to time $t=15$
is presented in Figure $\ref{F3}.$ Furthermore, in Table \ref{p3t1}, a
comparison of our results with those in the literature is given. One can see
from the table that present results are$\ $compatibly in good harmony with
their exact values and all of the compared ones.

\begin{table}[ptb]
\caption{Comparison of $\ $the calculated invariants of Problem 3 for
$h=k=0.1$ ($\mu=1,$ $c_{1}=4.5,$ $c_{2}=1.5,$ $c_{3}=0.5,$ $x_{1}=10,$
$x_{2}=25$, $x_{3}=35,$ $0\leq x\leq100,$ $0\leq t\leq15$ ).}%
\label{p3t1}
{\scriptsize
\begin{tabular}
[c]{cccccc}\hline
Method &  & $t$ & $I_{1}$ & $I_{2}$ & $I_{3}$\\\hline
CHCM-L &  & $0$ & \multicolumn{1}{l}{$77.999971$} &
\multicolumn{1}{l}{$655.277034$} & \multicolumn{1}{l}{$5451.148721$}\\
& \multicolumn{1}{l}{} & $3$ & \multicolumn{1}{l}{$78.000025$} &
\multicolumn{1}{l}{$651.326045$} & \multicolumn{1}{l}{$5384.366499$}\\
& \multicolumn{1}{l}{} & $6$ & \multicolumn{1}{l}{$78.000025$} &
\multicolumn{1}{l}{$655.118139$} & \multicolumn{1}{l}{$5449.115647$}\\
&  & $9$ & \multicolumn{1}{l}{$78.000025$} & \multicolumn{1}{l}{$655.286252$}
& \multicolumn{1}{l}{$5451.661801$}\\
&  & $12$ & \multicolumn{1}{l}{$78.000025$} & \multicolumn{1}{l}{$655.329978$}
& \multicolumn{1}{l}{$5451.907342$}\\
&  & $15$ & \multicolumn{1}{l}{$78.000020$} & \multicolumn{1}{l}{$655.337316$}
& \multicolumn{1}{l}{$5451.947083$}\\
CHCM-C &  & $15$ & \multicolumn{1}{l}{$77.999656$} &
\multicolumn{1}{l}{$655.329657$} & \multicolumn{1}{l}{$5451.857640$}\\
\cite{ek1} &  & $15$ & \multicolumn{1}{l}{$77.999994$} &
\multicolumn{1}{l}{$655.344625$} & \multicolumn{1}{l}{$5452.024410$}\\
\cite{bl3kraslan} &  & $15$ & \multicolumn{1}{l}{$78.00490$} &
\multicolumn{1}{l}{$652.3474$} & \multicolumn{1}{l}{$5412.232$}\\
\cite{bl3ali} & \multicolumn{1}{l}{$(h=0.4)$} & $15$ &
\multicolumn{1}{l}{$77.999984$} & \multicolumn{1}{l}{$652.411538$} &
\multicolumn{1}{l}{$5412.23185$}\\
\cite{bl3ydereli} & \multicolumn{1}{l}{$(h=0.5)$} & $15$ &
\multicolumn{1}{l}{$78.000222$} & \multicolumn{1}{l}{$655.341909$} &
\multicolumn{1}{l}{$5452.481409$}\\
\cite{bl3uddin} & $(h=0.1833,$ $k=0.05)$ & $15$ &
\multicolumn{1}{l}{$77.86967$} & \multicolumn{1}{l}{$654.09104$} &
\multicolumn{1}{l}{$5440.78956$}\\
\cite{bl3kaslan} & \multicolumn{1}{l}{} & $15$ &
\multicolumn{1}{l}{$77.995390$} & \multicolumn{1}{l}{$652.810400$} &
\multicolumn{1}{l}{$5411.6390$}\\\hline
Analytical & \multicolumn{1}{c}{} & \multicolumn{1}{c}{} &
\multicolumn{1}{l}{$78$} & \multicolumn{1}{l}{\ $655.2$} &
\multicolumn{1}{l}{$5450.4$}\\\hline
\end{tabular}
\ \ \ \ \ \ \ \ \ \ }\end{table}

\subsection{The Maxwellian initial condition}

\begin{figure}[ptb]
\begin{center}
\centering\includegraphics[width=0.48\textwidth]{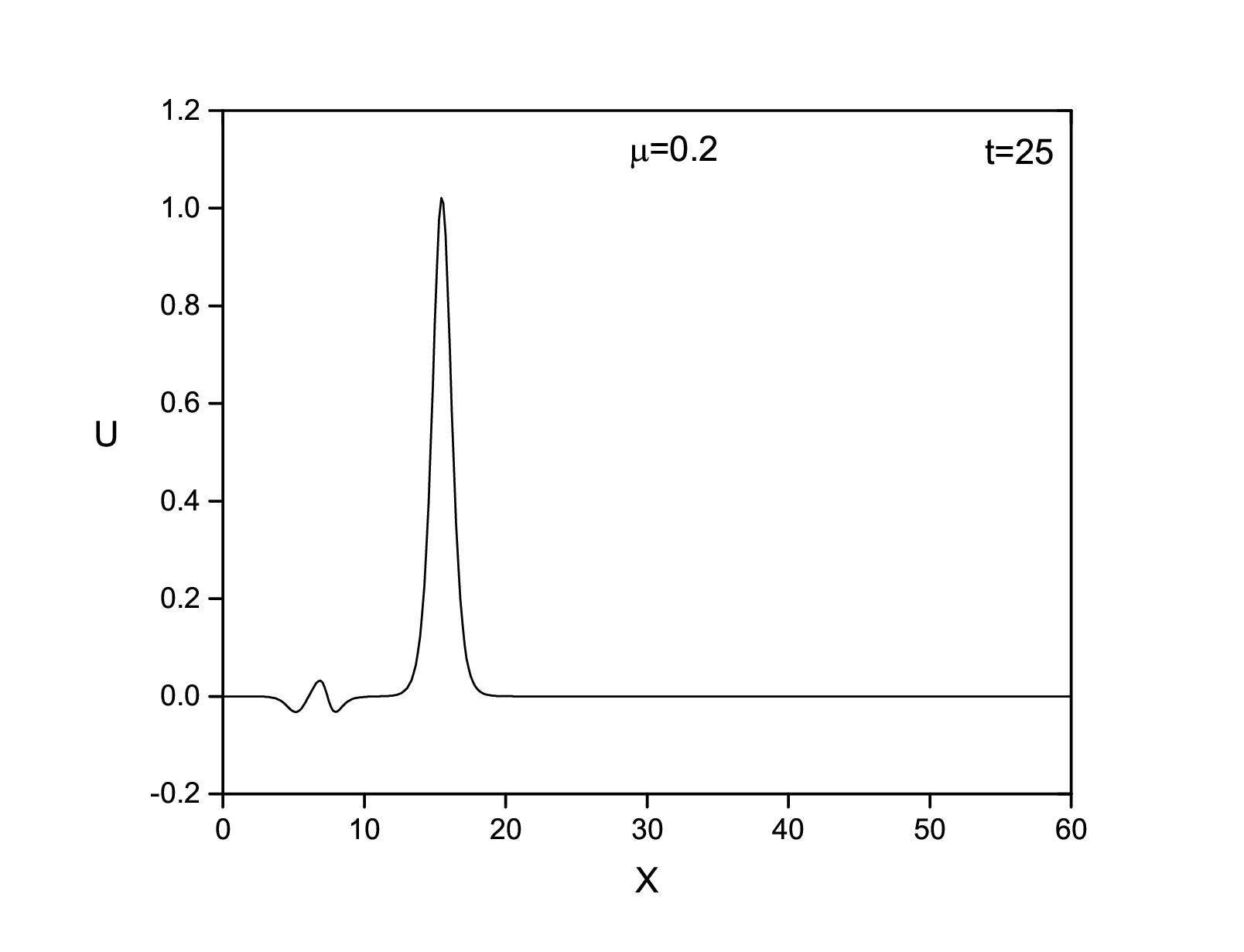}
\includegraphics[width=0.48\textwidth]{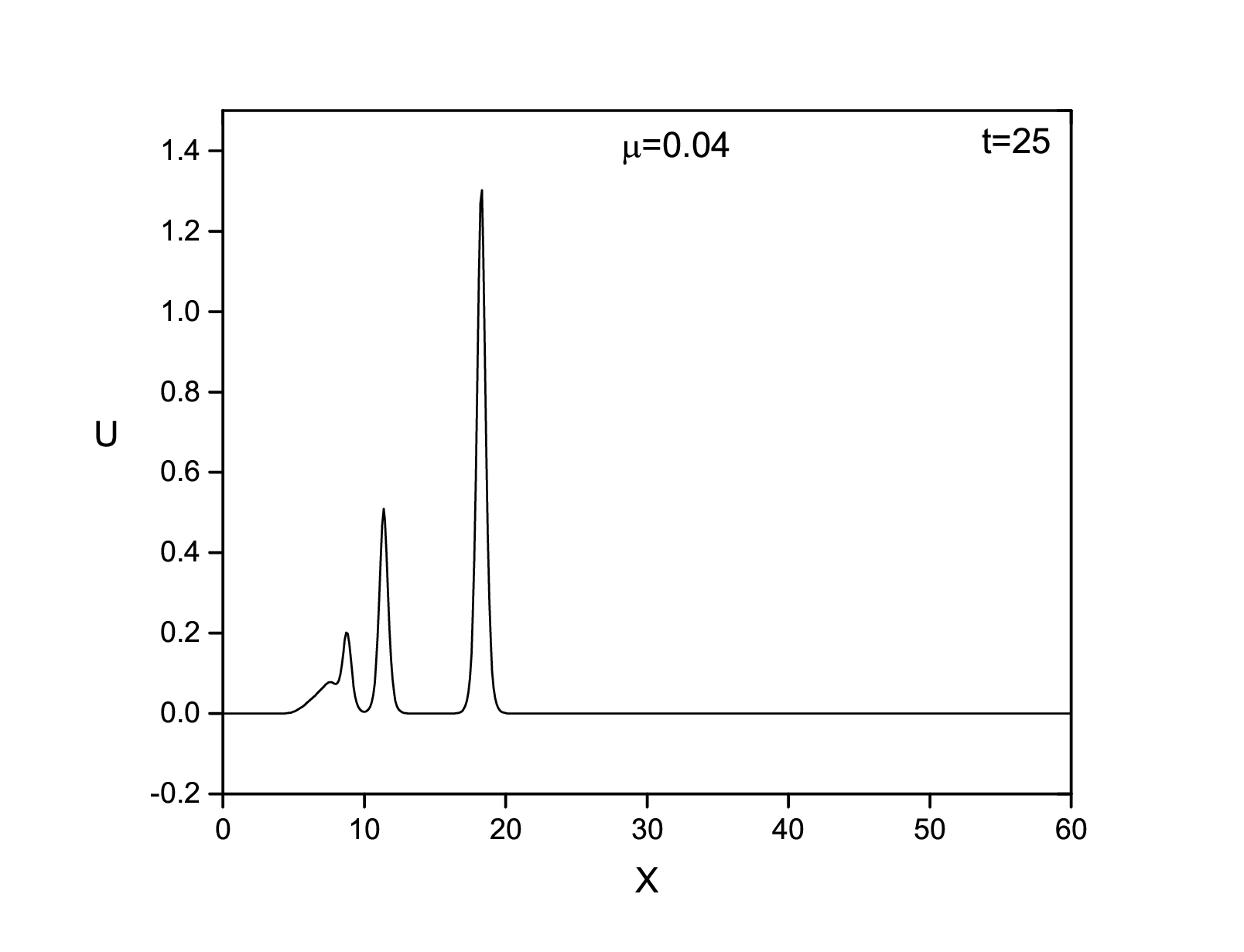}
\includegraphics[width=0.48\textwidth]{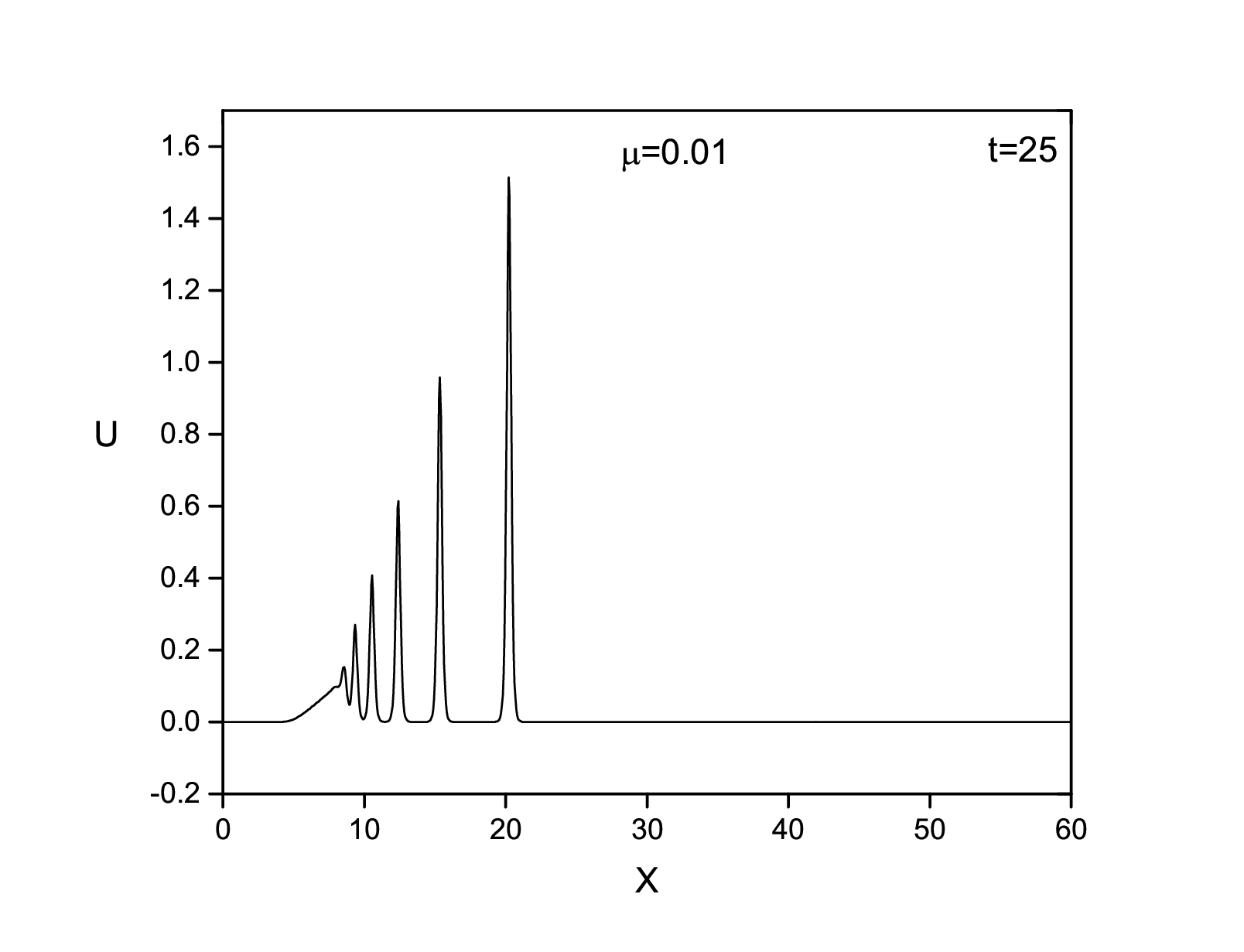}
\includegraphics[width=0.48\textwidth]{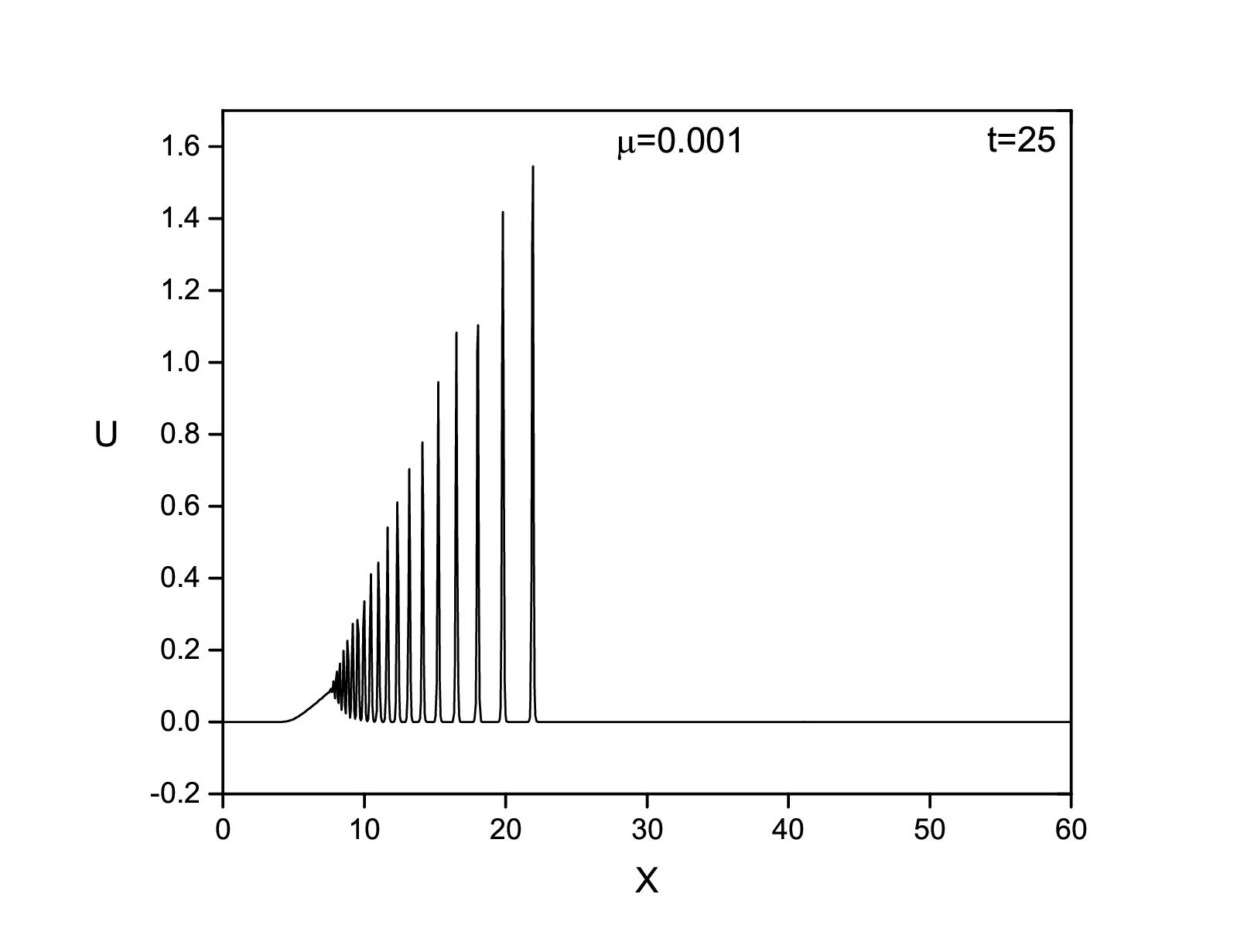}
\end{center}
\caption{The {\protect\footnotesize simula}%
\textrm{{\protect\footnotesize tions} {\protect\footnotesize for Maxwellian
initial condition }}{\protect\footnotesize for }${\protect\footnotesize \Delta
t=0.01}$.}%
\label{F4}%
\end{figure}

The fourth experimental problem dwells on the Maxwellian initial condition of
the following form\cite{bl3roshan}%
\[
U(x,0)=e^{\left(  -\left(  x-20\right)  ^{2}\right)  }.
\]
The simulations of the Maxwellian pulse are found for constant $\Delta t=0.01$
and different values of the $\mu=0.2,0.04,0.01$ and $0.001$, respectively$.$
Simulation of the waves for the values $\mu=0.2,$ $0.04,$\ $0.01$ and $0.001$
at $t=25$\ is presented in Figure $\ref{F4}.$ Moreover, in Table \ref{p4t1},
one can see a comparison of the present results with some of those given in
the literature. One can obviously see from the investigation of the table, the
newly obtained results are$\ $also in good harmony with the exact values and
all of the compared ones.

\begin{table}[ptb]
\caption{The computed invariants of Problem 4 and a comparion with those in
Refs.\cite{ek5} and \cite{bl3roshan} for values of $h=0.05$ and $k=0.025.$}%
\label{p4t1}
{\scriptsize
\begin{tabular}
[c]{ccccccccccccc}\hline
&  & \multicolumn{3}{c}{$I_{1}$} &  & \multicolumn{3}{c}{$I_{2}$} &  &
\multicolumn{3}{c}{$I_{3}$}\\\cline{3-5}\cline{7-9}\cline{11-13}%
$\mu$ & $t$ & CHCM-L & \cite{ek5} & \cite{bl3roshan} &  & CHCM-L & \cite{ek5}
& \cite{bl3roshan} &  & CHCM-L & \cite{ek5} & \cite{bl3roshan}\\\hline
\multicolumn{1}{l}{} & $0$ & \multicolumn{1}{l}{$1.77245$} &
\multicolumn{1}{l}{$1.77245$} & \multicolumn{1}{l}{$1.77245$} &
\multicolumn{1}{l}{} & \multicolumn{1}{l}{$1.37864$} &
\multicolumn{1}{l}{$1.37864$} & \multicolumn{1}{l}{$1.37864$} &
\multicolumn{1}{l}{} & \multicolumn{1}{l}{$1.02332$} &
\multicolumn{1}{l}{$1.02333$} & \multicolumn{1}{l}{$1.02333$}\\
\multicolumn{1}{l}{} & $3$ & \multicolumn{1}{l}{$1.77245$} &
\multicolumn{1}{l}{$1.77245$} & \multicolumn{1}{l}{$1.77245$} &
\multicolumn{1}{l}{} & \multicolumn{1}{l}{$1.37867$} &
\multicolumn{1}{l}{$1.37867$} & \multicolumn{1}{l}{$1.37923$} &
\multicolumn{1}{l}{} & \multicolumn{1}{l}{$1.02335$} &
\multicolumn{1}{l}{$1.02336$} & \multicolumn{1}{l}{$1.02355$}\\
$0.1$ & $6$ & \multicolumn{1}{l}{$1.77245$} & \multicolumn{1}{l}{$1.77245$} &
\multicolumn{1}{l}{$1.77245$} & \multicolumn{1}{l}{} &
\multicolumn{1}{l}{$1.37868$} & \multicolumn{1}{l}{$1.37868$} &
\multicolumn{1}{l}{$1.37880$} & \multicolumn{1}{l}{} &
\multicolumn{1}{l}{$1.02337$} & \multicolumn{1}{l}{$1.02337$} &
\multicolumn{1}{l}{$1.02338$}\\
& $9$ & \multicolumn{1}{l}{$1.77245$} & \multicolumn{1}{l}{$1.77245$} &
\multicolumn{1}{l}{$1.77245$} & \multicolumn{1}{l}{} &
\multicolumn{1}{l}{$1.37868$} & \multicolumn{1}{l}{$1.37869$} &
\multicolumn{1}{l}{$1.37877$} & \multicolumn{1}{l}{} &
\multicolumn{1}{l}{$1.02337$} & \multicolumn{1}{l}{$1.02338$} &
\multicolumn{1}{l}{$1.02336$}\\
\multicolumn{1}{l}{} & $12$ & \multicolumn{1}{l}{$1.77245$} &
\multicolumn{1}{l}{$1.77245$} & \multicolumn{1}{l}{$1.77245$} &
\multicolumn{1}{l}{} & \multicolumn{1}{l}{$1.37868$} &
\multicolumn{1}{l}{$1.37869$} & \multicolumn{1}{l}{$1.37885$} &
\multicolumn{1}{l}{} & \multicolumn{1}{l}{$1.02337$} &
\multicolumn{1}{l}{$1.02338$} & \multicolumn{1}{l}{$1.02339$}\\
&  & \multicolumn{1}{l}{} & \multicolumn{1}{l}{} & \multicolumn{1}{l}{} &
\multicolumn{1}{l}{} & \multicolumn{1}{l}{} & \multicolumn{1}{l}{} &
\multicolumn{1}{l}{} & \multicolumn{1}{l}{} & \multicolumn{1}{l}{} &
\multicolumn{1}{l}{} & \multicolumn{1}{l}{}\\
& $0$ & \multicolumn{1}{l}{$1.77245$} & \multicolumn{1}{l}{$1.77245$} &
\multicolumn{1}{l}{$1.77245$} & \multicolumn{1}{l}{} &
\multicolumn{1}{l}{$1.31597$} & \multicolumn{1}{l}{$1.31598$} &
\multicolumn{1}{l}{$1.31598$} & \multicolumn{1}{l}{} &
\multicolumn{1}{l}{$1.02332$} & \multicolumn{1}{l}{$1.02333$} &
\multicolumn{1}{l}{$1.02333$}\\
& $3$ & \multicolumn{1}{l}{$1.77245$} & \multicolumn{1}{l}{$1.77245$} &
\multicolumn{1}{l}{$1.77245$} & \multicolumn{1}{l}{} &
\multicolumn{1}{l}{$1.31606$} & \multicolumn{1}{l}{$1.31606$} &
\multicolumn{1}{l}{$1.31648$} & \multicolumn{1}{l}{} &
\multicolumn{1}{l}{$1.02345$} & \multicolumn{1}{l}{$1.02345$} &
\multicolumn{1}{l}{$1.02356$}\\
$0.05$ & $6$ & \multicolumn{1}{l}{$1.77245$} & \multicolumn{1}{l}{$1.77245$} &
\multicolumn{1}{l}{$1.77245$} & \multicolumn{1}{l}{} &
\multicolumn{1}{l}{$1.31611$} & \multicolumn{1}{l}{$1.31611$} &
\multicolumn{1}{l}{$1.31619$} & \multicolumn{1}{l}{} &
\multicolumn{1}{l}{$1.02352$} & \multicolumn{1}{l}{$1.02352$} &
\multicolumn{1}{l}{$1.02340$}\\
& $9$ & \multicolumn{1}{l}{$1.77245$} & \multicolumn{1}{l}{$1.77245$} &
\multicolumn{1}{l}{$1.77245$} & \multicolumn{1}{l}{} &
\multicolumn{1}{l}{$1.31612$} & \multicolumn{1}{l}{$1.31611$} &
\multicolumn{1}{l}{$1.31617$} & \multicolumn{1}{l}{} &
\multicolumn{1}{l}{$1.02353$} & \multicolumn{1}{l}{$1.02353$} &
\multicolumn{1}{l}{$1.02339$}\\
& $12$ & \multicolumn{1}{l}{$1.77245$} & \multicolumn{1}{l}{$1.77245$} &
\multicolumn{1}{l}{$1.77245$} & \multicolumn{1}{l}{} &
\multicolumn{1}{l}{$1.31612$} & \multicolumn{1}{l}{$1.31611$} &
\multicolumn{1}{l}{$1.31612$} & \multicolumn{1}{l}{} &
\multicolumn{1}{l}{$1.02353$} & \multicolumn{1}{l}{$1.02353$} &
\multicolumn{1}{l}{$1.02339$}\\
&  & \multicolumn{1}{l}{} & \multicolumn{1}{l}{} & \multicolumn{1}{l}{} &
\multicolumn{1}{l}{} & \multicolumn{1}{l}{} & \multicolumn{1}{l}{} &
\multicolumn{1}{l}{} & \multicolumn{1}{l}{} & \multicolumn{1}{l}{} &
\multicolumn{1}{l}{} & \multicolumn{1}{l}{}\\
& $0$ & \multicolumn{1}{l}{$1.77245$} & \multicolumn{1}{l}{$1.77245$} &
\multicolumn{1}{l}{$1.77245$} & \multicolumn{1}{l}{} &
\multicolumn{1}{l}{$1.28464$} & \multicolumn{1}{l}{$1.28464$} &
\multicolumn{1}{l}{$1.28464$} & \multicolumn{1}{l}{} &
\multicolumn{1}{l}{$1.02332$} & \multicolumn{1}{l}{$1.02333$} &
\multicolumn{1}{l}{$1.02333$}\\
& $3$ & \multicolumn{1}{l}{$1.77245$} & \multicolumn{1}{l}{$1.77238$} &
\multicolumn{1}{l}{$1.77245$} & \multicolumn{1}{l}{} &
\multicolumn{1}{l}{$1.28488$} & \multicolumn{1}{l}{$1.28487$} &
\multicolumn{1}{l}{$1.28520$} & \multicolumn{1}{l}{} &
\multicolumn{1}{l}{$1.02373$} & \multicolumn{1}{l}{$1.02372$} &
\multicolumn{1}{l}{$1.02357$}\\
$0.025$ & $6$ & \multicolumn{1}{l}{$1.77245$} & \multicolumn{1}{l}{$1.77254$}
& \multicolumn{1}{l}{$1.77245$} & \multicolumn{1}{l}{} &
\multicolumn{1}{l}{$1.28500$} & \multicolumn{1}{l}{$1.28499$} &
\multicolumn{1}{l}{$1.28492$} & \multicolumn{1}{l}{} &
\multicolumn{1}{l}{$1.02390$} & \multicolumn{1}{l}{$1.02392$} &
\multicolumn{1}{l}{$1.02340$}\\
& $9$ & \multicolumn{1}{l}{$1.77245$} & \multicolumn{1}{l}{$1.77233$} &
\multicolumn{1}{l}{$1.77245$} & \multicolumn{1}{l}{} &
\multicolumn{1}{l}{$1.28501$} & \multicolumn{1}{l}{$1.28496$} &
\multicolumn{1}{l}{$1.28418$} & \multicolumn{1}{l}{} &
\multicolumn{1}{l}{$1.02391$} & \multicolumn{1}{l}{$1.02386$} &
\multicolumn{1}{l}{$1.02337$}\\
& $12$ & \multicolumn{1}{l}{$1.77245$} & \multicolumn{1}{l}{$1.77253$} &
\multicolumn{1}{l}{$1.77245$} & \multicolumn{1}{l}{} &
\multicolumn{1}{l}{$1.28501$} & \multicolumn{1}{l}{$1.28497$} &
\multicolumn{1}{l}{$1.28474$} & \multicolumn{1}{l}{} &
\multicolumn{1}{l}{$1.02391$} & \multicolumn{1}{l}{$1.02390$} &
\multicolumn{1}{l}{$1.02337$}\\
&  & \multicolumn{1}{l}{} & \multicolumn{1}{l}{} & \multicolumn{1}{l}{} &
\multicolumn{1}{l}{} & \multicolumn{1}{l}{} & \multicolumn{1}{l}{} &
\multicolumn{1}{l}{} & \multicolumn{1}{l}{} & \multicolumn{1}{l}{} &
\multicolumn{1}{l}{} & \multicolumn{1}{l}{}\\
& $0$ & \multicolumn{1}{l}{$1.77245$} & \multicolumn{1}{l}{$1.77245$} &
\multicolumn{1}{l}{$1.77245$} & \multicolumn{1}{l}{} &
\multicolumn{1}{l}{$1.26584$} & \multicolumn{1}{l}{$1.26585$} &
\multicolumn{1}{l}{$1.26585$} & \multicolumn{1}{l}{} &
\multicolumn{1}{l}{$1.02332$} & \multicolumn{1}{l}{$1.02333$} &
\multicolumn{1}{l}{$1.02333$}\\
& $3$ & \multicolumn{1}{l}{$1.77245$} & \multicolumn{1}{l}{$1.77247$} &
\multicolumn{1}{l}{$1.77245$} & \multicolumn{1}{l}{} &
\multicolumn{1}{l}{$1.26666$} & \multicolumn{1}{l}{$1.26572$} &
\multicolumn{1}{l}{$1.26632$} & \multicolumn{1}{l}{} &
\multicolumn{1}{l}{$1.02481$} & \multicolumn{1}{l}{$1.02293$} &
\multicolumn{1}{l}{$1.02330$}\\
$0.01$ & $6$ & \multicolumn{1}{l}{$1.77245$} & \multicolumn{1}{l}{$1.77253$} &
\multicolumn{1}{l}{$1.77245$} & \multicolumn{1}{l}{} &
\multicolumn{1}{l}{$1.26695$} & \multicolumn{1}{l}{$1.26579$} &
\multicolumn{1}{l}{$1.26599$} & \multicolumn{1}{l}{} &
\multicolumn{1}{l}{$1.02520$} & \multicolumn{1}{l}{$1.02297$} &
\multicolumn{1}{l}{$1.02294$}\\
& $9$ & \multicolumn{1}{l}{$1.77245$} & \multicolumn{1}{l}{$1.77252$} &
\multicolumn{1}{l}{$1.77245$} & \multicolumn{1}{l}{} &
\multicolumn{1}{l}{$1.26700$} & \multicolumn{1}{l}{$1.26562$} &
\multicolumn{1}{l}{$1.26639$} & \multicolumn{1}{l}{} &
\multicolumn{1}{l}{$1.02523$} & \multicolumn{1}{l}{$1.02295$} &
\multicolumn{1}{l}{$1.02295$}\\
& $12$ & \multicolumn{1}{l}{$1.77245$} & \multicolumn{1}{l}{$1.77253$} &
\multicolumn{1}{l}{$1.77245$} & \multicolumn{1}{l}{} &
\multicolumn{1}{l}{$1.26702$} & \multicolumn{1}{l}{$1.26566$} &
\multicolumn{1}{l}{$1.26567$} & \multicolumn{1}{l}{} &
\multicolumn{1}{l}{$1.02523$} & \multicolumn{1}{l}{$1.02292$} &
\multicolumn{1}{l}{$1.02293$}\\\hline
\end{tabular}
}\end{table}

\subsection{Undular Bore}

In the fifth experimental problem, the EW equation (\ref{1}) is taken into
consideration in the finite range $a\leq x\leq b$ with the boundary conditions%
\begin{align*}
U(a,t)  &  =U_{0}\\
U(b,t)  &  =0
\end{align*}
and the initial condition
\[
U(x,0)=0.5U_{0}\left[  1-\tanh(\frac{x-x_{0}}{d})\right]
\]
to examine undular bore formation \cite{bl3bsaka}. In this equation
$U(x,0)\ $stands for the height of the water on the stagnant water at initial
time, $d$ stands for the difference in slopes between the deep and stagnant
water. The change in the water level $U(x,0)$ occurrs at the point $x=x_{0}$.
The stagnant water can be obesreved to the right hand of the zone and at the
additional elevation $U_{0}$ from the surface $U=0$ the flow of water moves
from the left into the stagnant water.

In this experimental problem, the conservation constants of $I_{1},I_{2}$ and
$I_{3}$ do not remain constant however linearly increase in the following
ratios $M_{1},M_{2}$ and $M_{3}$, respectively \cite{bl3gardnerayoup}.%
\begin{align*}
M_{1}  &  =\frac{d}{dt}I_{1}=\frac{d}{dt}\int_{a}^{b}Udx=\frac{1}{2}%
(U_{0})^{2},\\
M_{2}  &  =\frac{d}{dt}I_{2}=\frac{d}{dt}\int_{a}^{b}[U^{2}+\mu(U_{x}%
)^{2}]dx=\frac{2}{3}(U_{0})^{3},\\
M_{3}  &  =\frac{d}{dt}I_{3}=\frac{d}{dt}\int_{a}^{b}U^{3}dx=\frac{3}{4}%
(U_{0})^{4}.
\end{align*}
During numerical computations the values $U_{0}=0.1,$ $\mu=0.16666667$ and
$x_{0}=0$ are utilized. Therefore, the linearly increasing ratios of the
conservation constants for those parameters are found as%
\[
M_{1}=5e-3,\text{ }M_{2}=6.66667e-4,\text{ }M_{3}=7.5e-5.
\]
\begin{table}[ptb]
\caption{Comparison of $\ $the computed invariants of Problem 5 for
$k=0.05$\ and $h=0.07$ ($\mu=0.16666667,$ $d=2,$ $x_{0}=0$, $U_{0}=0.1,$
$0\leq t\leq800,$ $-20\leq x\leq50$)}%
\label{p5t1}%
{\scriptsize
\begin{tabular}
[c]{ccccccc}\hline
Method & $t$ & $I_{1}$ & $I_{2}$ & $I_{3}$ & $x$ & $U$\\\hline
CHCM-L & $0$ & \multicolumn{1}{l}{$1.996500$} & \multicolumn{1}{l}{$0.189927$}
& \multicolumn{1}{l}{$0.018465$} & \multicolumn{1}{l}{$-20.00$} &
\multicolumn{1}{l}{$0.10000$}\\
\multicolumn{1}{l}{} & $100$ & \multicolumn{1}{l}{$2.496499$} &
\multicolumn{1}{l}{$0.256594$} & \multicolumn{1}{l}{$0.025965$} &
\multicolumn{1}{l}{$3.73$} & \multicolumn{1}{l}{$0.15730$}\\
& $200$ & \multicolumn{1}{l}{$2.996499$} & \multicolumn{1}{l}{$0.323261$} &
\multicolumn{1}{l}{$0.033465$} & \multicolumn{1}{l}{$9.40$} &
\multicolumn{1}{l}{$0.17606$}\\
\multicolumn{1}{l}{} & $300$ & \multicolumn{1}{l}{$3.496499$} &
\multicolumn{1}{l}{$0.389928$} & \multicolumn{1}{l}{$0.040965$} &
\multicolumn{1}{l}{$15.35$} & \multicolumn{1}{l}{$0.18010$}\\
& $400$ & \multicolumn{1}{l}{$3.996499$} & \multicolumn{1}{l}{$0.456595$} &
\multicolumn{1}{l}{$0.048465$} & \multicolumn{1}{l}{$21.37$} &
\multicolumn{1}{l}{$0.18214$}\\
& $500$ & \multicolumn{1}{l}{$4.496499$} & \multicolumn{1}{l}{$0.523262$} &
\multicolumn{1}{l}{$0.055965$} & \multicolumn{1}{l}{$27.46$} &
\multicolumn{1}{l}{$0.18321$}\\
& $600$ & \multicolumn{1}{l}{$4.996499$} & \multicolumn{1}{l}{$0.589929$} &
\multicolumn{1}{l}{$0.063465$} & \multicolumn{1}{l}{$33.55$} &
\multicolumn{1}{l}{$0.18378$}\\
& $700$ & \multicolumn{1}{l}{$5.496499$} & \multicolumn{1}{l}{$0.656596$} &
\multicolumn{1}{l}{$0.070965$} & \multicolumn{1}{l}{$39.71$} &
\multicolumn{1}{l}{$0.18440$}\\
& $800$ & \multicolumn{1}{l}{$5.996475$} & \multicolumn{1}{l}{$0.723263$} &
\multicolumn{1}{l}{$0.078465$} & \multicolumn{1}{l}{$45.87$} &
\multicolumn{1}{l}{$0.18474$}\\
CHCM-C & $800$ & \multicolumn{1}{l}{$5.995270$} &
\multicolumn{1}{l}{$0.722959$} & \multicolumn{1}{l}{$0.078422$} &
\multicolumn{1}{l}{$45.87$} & \multicolumn{1}{l}{$0.18467$}\\
\cite{ek1} & $800$ & \multicolumn{1}{l}{$6.003322$} &
\multicolumn{1}{l}{$0.723860$} & \multicolumn{1}{l}{$0.078533$} &
\multicolumn{1}{l}{$45.87$} & \multicolumn{1}{l}{$0.18451$}\\
\cite{bl3gardnerayoup} & $800$ & \multicolumn{1}{l}{$5.994366$} &
\multicolumn{1}{l}{$0.712677$} & \multicolumn{1}{l}{$0.076876$} &
\multicolumn{1}{l}{$45.70$} & \multicolumn{1}{l}{$0.183918$}\\
\cite{bl014} & $800$ & \multicolumn{1}{l}{$5.996473$} &
\multicolumn{1}{l}{$0.722126$} & \multicolumn{1}{l}{$0.078465$} &
\multicolumn{1}{l}{$45.87$} & \multicolumn{1}{l}{$0.184431$}\\
\cite{bl3esen} & $800$ & \multicolumn{1}{l}{$6.003478$} &
\multicolumn{1}{l}{$0.723605$} & \multicolumn{1}{l}{$0.078426$} &
\multicolumn{1}{l}{$45.87$} & \multicolumn{1}{l}{$0.184518$}\\
\cite{bl1idrisdag} & $800$ & \multicolumn{1}{l}{$6.003194$} &
\multicolumn{1}{l}{$0.723867$} & \multicolumn{1}{l}{$0.078534$} &
\multicolumn{1}{l}{$45.85$} & \multicolumn{1}{l}{$0.18460$}\\
\cite{bl3dogan} & $800$ & \multicolumn{1}{l}{$5.669824$} &
\multicolumn{1}{l}{$0.660997$} & \multicolumn{1}{l}{$0.070677$} &
\multicolumn{1}{l}{$46.73$} & \multicolumn{1}{l}{$0.197568$}\\
\cite{bl3bsaka} & $800$ & \multicolumn{1}{l}{$6.00248$} &
\multicolumn{1}{l}{$0.72402$} & \multicolumn{1}{l}{$0.07853$} &
\multicolumn{1}{l}{$45.85$} & \multicolumn{1}{l}{$0.184713$}\\
\cite{bl3bsakadag} & $800$ & \multicolumn{1}{l}{$6.002474$} &
\multicolumn{1}{l}{$0.723860$} & \multicolumn{1}{l}{$0.078525$} &
\multicolumn{1}{l}{$45.85$} & \multicolumn{1}{l}{$0.18471$}\\\hline
\end{tabular}
\ \ \ \ \ \ \ \ \ \ }\end{table}

\begin{figure}[ptb]
\begin{center}
\centering\includegraphics[width=0.48\textwidth]{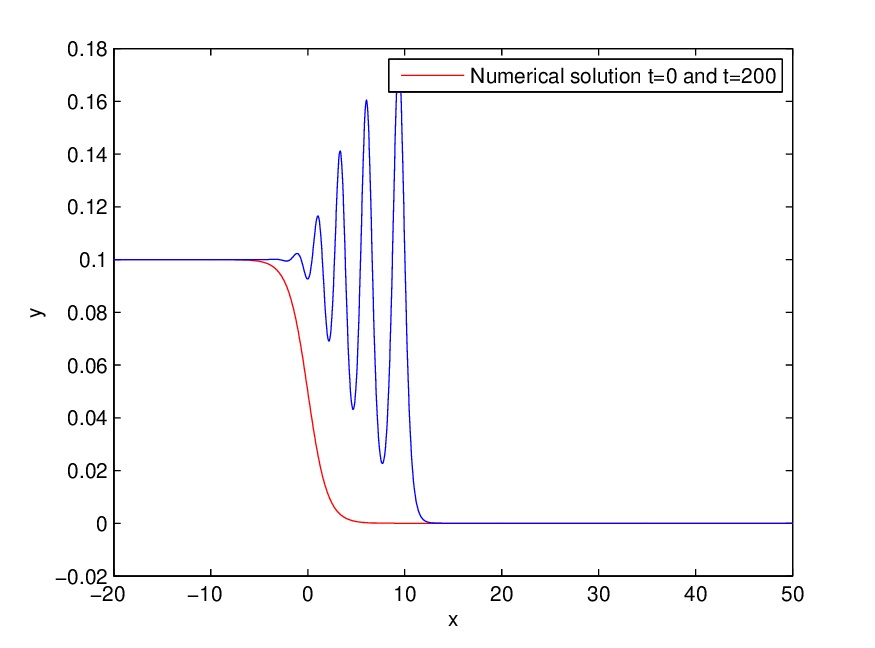}
\includegraphics[width=0.48\textwidth]{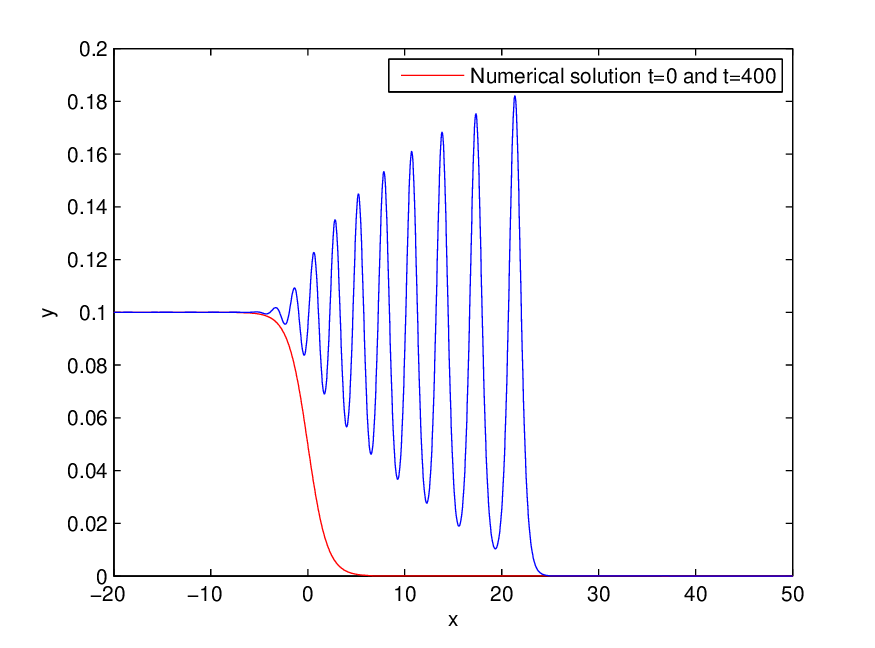}
\includegraphics[width=0.48\textwidth]{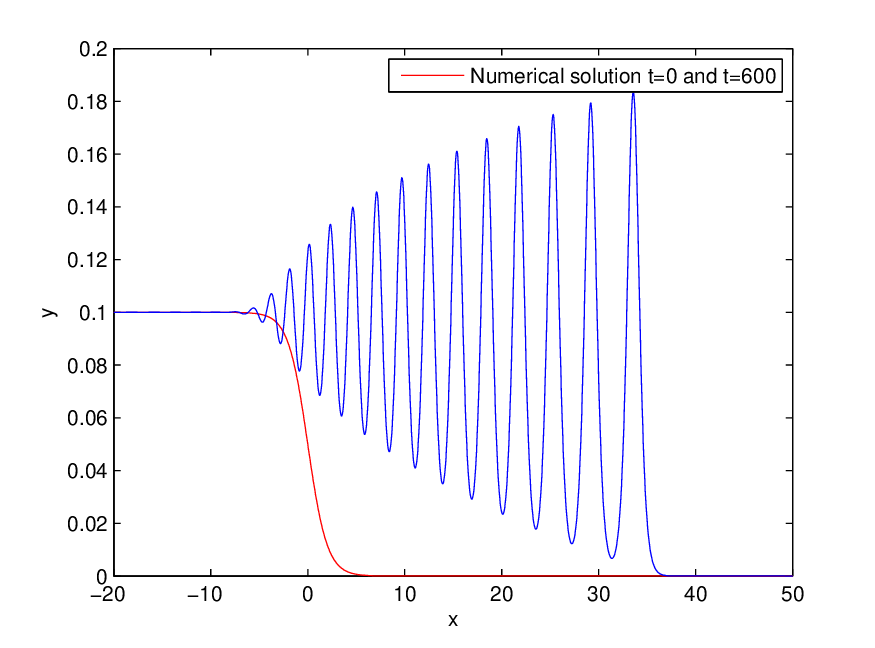}
\includegraphics[width=0.48\textwidth]{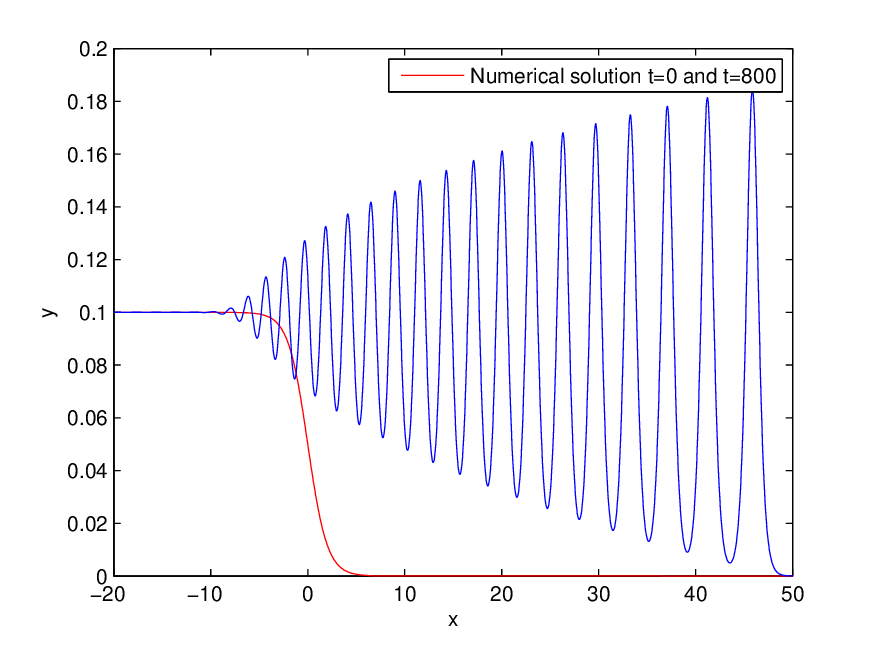}
\end{center}
\caption{The profiles {\protect\footnotesize and undulation profiles for
}${\protect\footnotesize d=2}$ at different times.}%
\label{F5}%
\end{figure}

The simulation process for the undular bore at different times $t$ and
$d=2$\ is presented in Figure $\ref{F5}.$ Furthermore, in Table \ref{p5t1},
the present results have been compared to some of those available in the
literature. One can easily see in this table that the presented method
produces good results and they are also are$\ $in very good harmony with both
their exact values and all of the compared ones.

\subsection{Soliton collision}

In the sixth experimental problem, the collision of two solitary waves having
the initial condition \cite{bl013}%

\[
U(x,0)=\overset{2}{\underset{j=1}{\sum}}3c_{j}\sec h^{2}\left[  \frac{1}%
{2}\left(  x-x_{j}-c_{j}\right)  \right]  .
\]
will be considered.

\begin{figure}[ptb]
\begin{center}
\centering\includegraphics[width=0.48\textwidth]{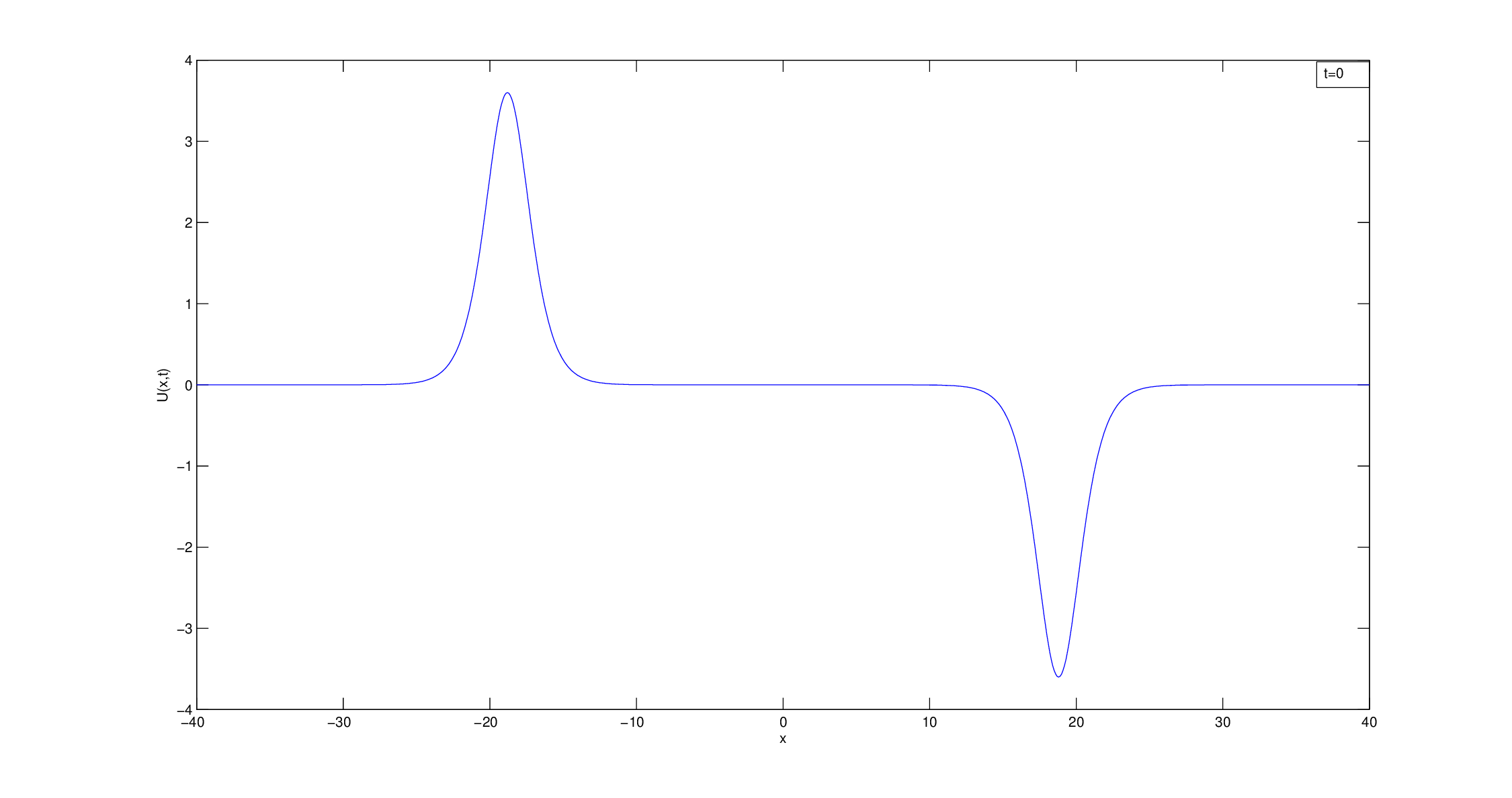}
\includegraphics[width=0.48\textwidth]{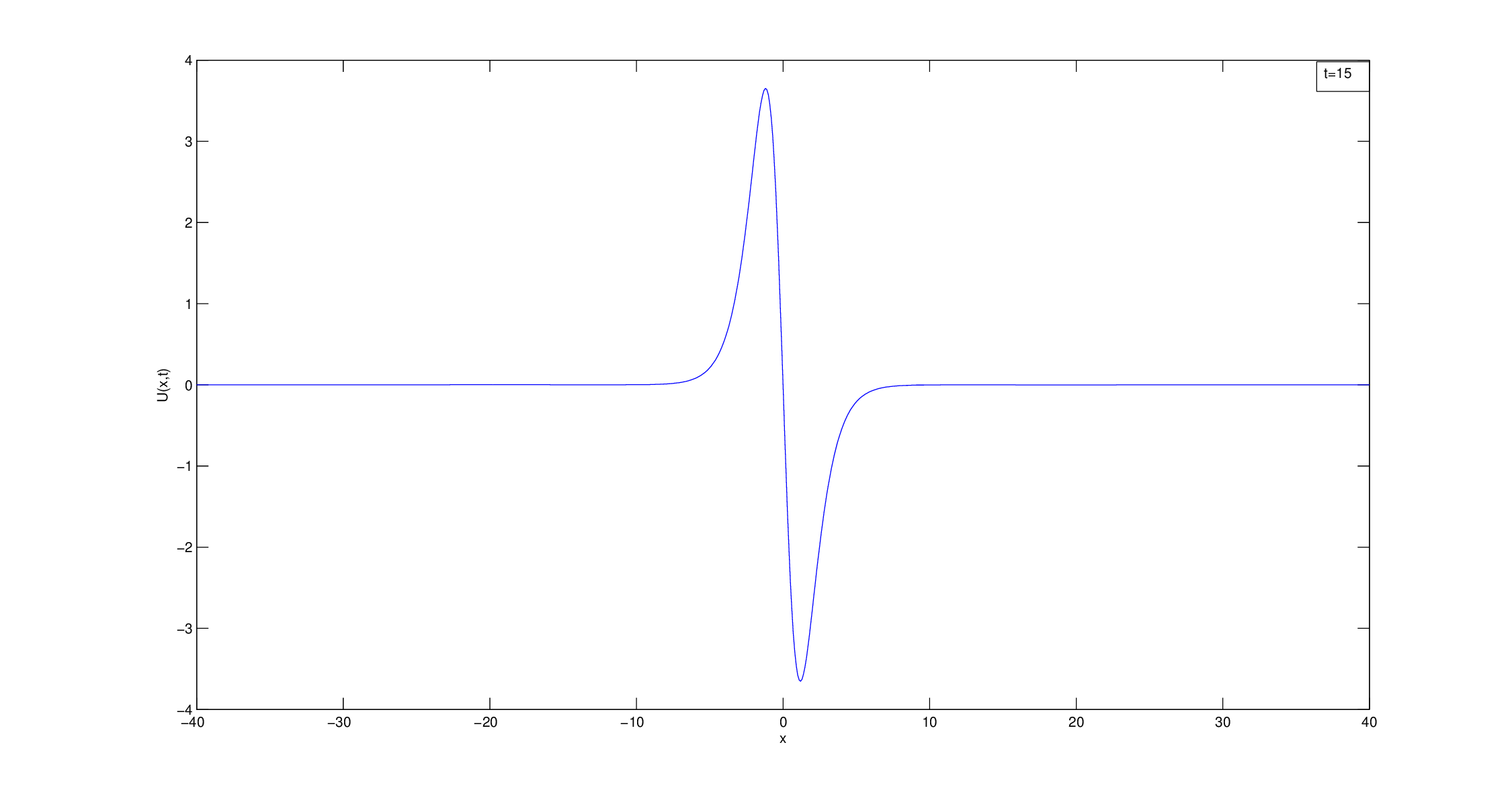}
\includegraphics[width=0.48\textwidth]{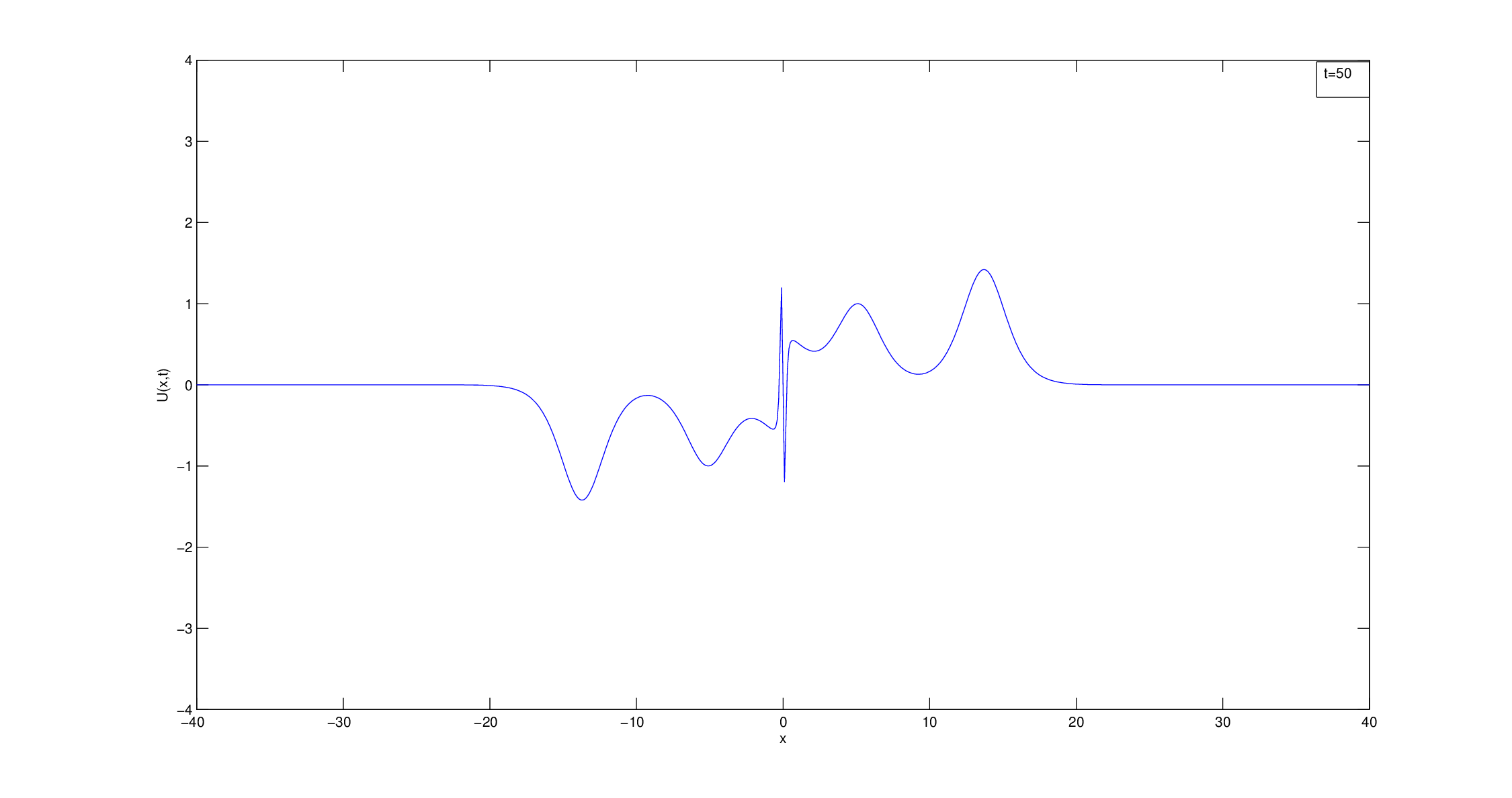}
\includegraphics[width=0.48\textwidth]{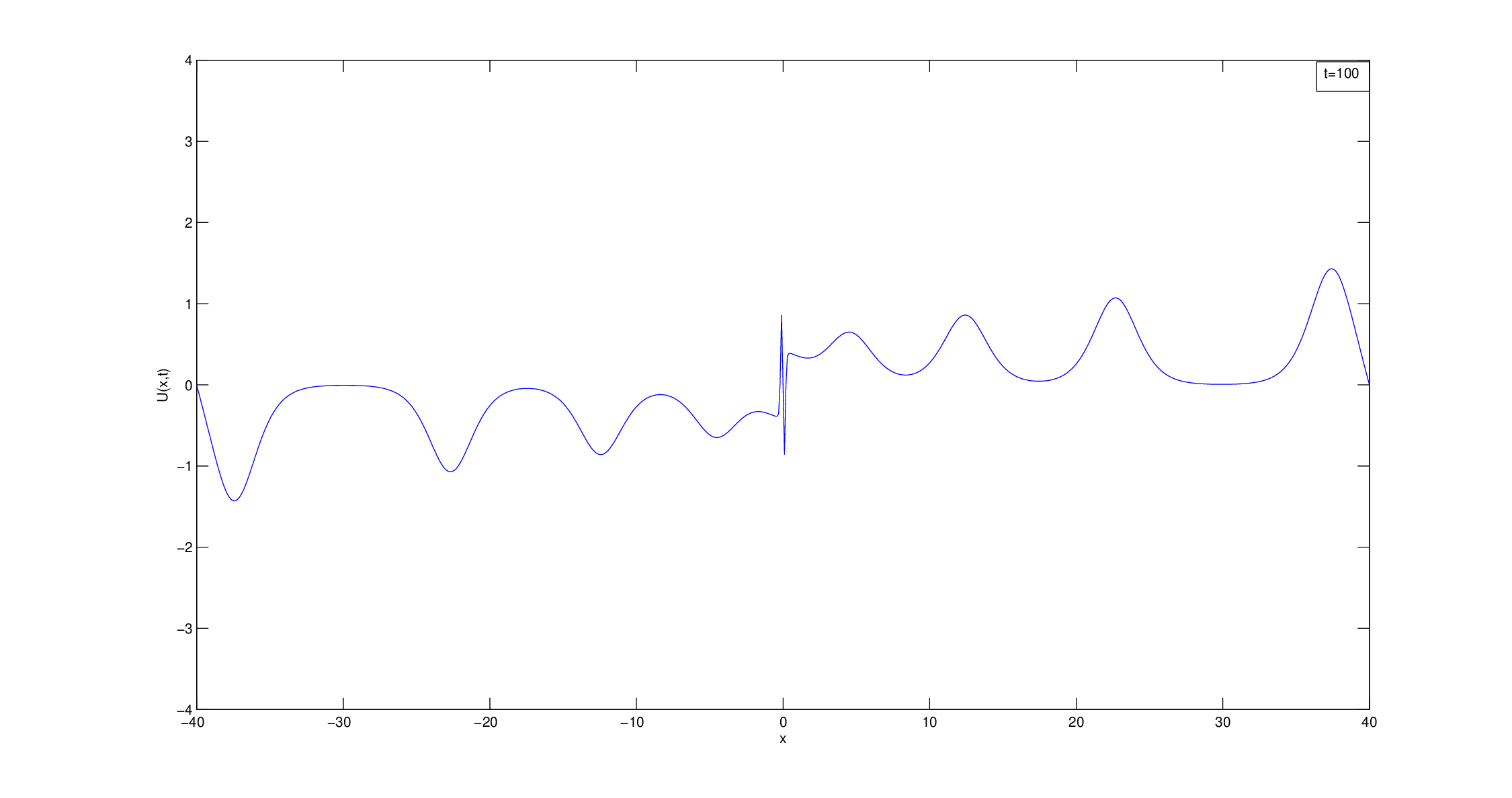}
\end{center}
\caption{{\protect\footnotesize Clash of two solitary waves.}}%
\label{F6}%
\end{figure}

These solitary waves are also presented like in the phenomena of interaction
of $2$ solitary waves except the fact that their signs are different and move
toward to one another. At collision time, a singularity happens and leaves
smaller waves behind However, when time elapses, these small singularities die
out. For the sake of computational aims, the following parameters $c_{1}%
=-1.2$, $c_{2}=1.2$, $\mu=1,$ $x_{1}=-20,$ $x_{2}=20$ with $\Delta t=0.1$ are
used over the solution domain $\left[  -40,40\right]  .$ The simulation
process for the collision of solitons for different values of $t=0,15,50,100$
is illustrated in Figure $\ref{F6}.$ One can see from this figure that the
waves display the expected physical behavior of the problem.

\section{\textbf{Conclusion }}

The numerical solutions of the EW equation which can be seen as an alternative
to the well-known KdV equation are found using cubic Hermite B-spline
collocation finite element method. To be able to establish the efficiency and
accuracy of the presented method with the help of the Crank-Nicolson type
approximation its validity, six test problems are considered and the obtained
results are tested by comparing with the previosly published ones especially
using the error norms $L_{2}$ and $L_{\infty}$. It is seen from all the
computed results that the presented numerical scheme produces reasonable
accurate results which are also in good agreement with exact ones and also
those of other researchers for the same parameters. As a future work, the
currently presented method may also easily and successfully be used to find
the numerical solutions of other frequenltly used non-linear PDEs seen in
varioud branches of mathematics and science that have a crucial role in
modelling natural phenomena.

\subsection*{Author contributions}

All persons who meet authorship criteria are listed as authors, and all
authors certify that they have participated sufficiently in the work to take
public responsibility for the content, including participation in the concept,
design, analysis, writing, or revision of the manuscript. Furthermore, each
author certifies that this material or similar material has not been and will
not be submitted to or published in any other publication.

\subsection*{Financial disclosure}

There are no funders to report for this submission.

\subsection*{Conflict of interest}

The authors declare that there is no conflict of interests regarding the
publication of this article.

\end{document}